\documentclass[11pt,a4paper]{article}

\usepackage[margin=1in]{geometry}
\usepackage{palatino}
\usepackage{subcaption}

\usepackage{enumitem}

\usepackage[
  colorlinks,
  linkcolor = blue,
  citecolor = blue,
  urlcolor = blue]{hyperref}

\usepackage{amsmath,amssymb,amsthm}
\usepackage{mathtools}
\mathtoolsset{centercolon}

\usepackage[affil-it]{authblk}
\usepackage{tabu}

\newcommand{\N}{\mathbb{N}}

\newcommand{\x}{\otimes}

\newcommand{\be}{\begin{equation}}
\newcommand{\ee}{\end{equation}}

\newcommand{\Thm}[1]{\hyperref[thm:#1]{Theorem~\ref*{thm:#1}}}
\newcommand{\Lem}[1]{\hyperref[lem:#1]{Lemma~\ref*{lem:#1}}}
\newcommand{\Cor}[1]{\hyperref[cor:#1]{Corollary~\ref*{cor:#1}}}
\newcommand{\Def}[1]{\hyperref[def:#1]{Definition~\ref*{def:#1}}}
\newcommand{\Obs}[1]{\hyperref[obs:#1]{Observation~\ref*{obs:#1}}}
\newcommand{\Prop}[1]{\hyperref[prop:#1]{Proposition~\ref*{prop:#1}}}
\newcommand{\Rem}[1]{\hyperref[rem:#1]{Remark~\ref*{rem:#1}}}
\newcommand{\Ex}[1]{\hyperref[ex:#1]{Example~\ref*{ex:#1}}}

\newcommand{\Sec}[1]{\hyperref[sec:#1]{Section~\ref*{sec:#1}}}

\newcommand{\Fig}[1]{\hyperref[fig:#1]{Figure~\ref*{fig:#1}}}
\newcommand{\Tab}[1]{\hyperref[tab:#1]{Table~\ref*{tab:#1}}}
\newcommand{\EqRef}[1]{\hyperref[eq:#1]{(\ref*{eq:#1})}}
\newcommand{\Eq}[1]{Equation~\hyperref[eq:#1]{(\ref*{eq:#1})}}

\makeatletter
\newtheorem*{rep@theorem}{\rep@title}
\newcommand{\newreptheorem}[2]{%
\newenvironment{rep#1}[1]{%
 \def\rep@title{#2 \ref{##1}}%
 \begin{rep@theorem}}%
 {\end{rep@theorem}}}
\makeatother

\newreptheorem{theorem}{Theorem}
\newreptheorem{lemma}{Lemma}

\newtheorem{theorem}{Theorem}[section]
\newtheorem*{theorem*}{Theorem}
\newtheorem{lemma}[theorem]{Lemma}
\newtheorem{cor}[theorem]{Corollary}

\newtheorem{prop}[theorem]{Proposition}

\theoremstyle{definition}
\newtheorem{definition}[theorem]{Definition}
\newtheorem{remark}[theorem]{Remark}
\newtheorem{example}[theorem]{Example}

\DeclareMathOperator{\qut}{Qut}

\newcommand{\qeds}{\qed\vspace{.2cm}}

\newcommand{\p}{\ensuremath{\mathbf{p}}}

\newcommand{\E}{\ensuremath{\mathcal{E}}}

\usepackage{blkarray}
\usepackage{comment}
\usepackage{tikz}

\title{Solution group representations as quantum symmetries of graphs}%Change the title?

\author[1]{David E.~Roberson\footnote{davideroberson@gmail.com}} 
\author[2]{Simon Schmidt\footnote{sisc@math.ku.dk}}

\affil[1]{Department of Applied Mathematics and Computer Science, Technical University of Denmark, DK-2800 Lyngby, Denmark}
\affil[2]{QMATH, Department of Mathematical Sciences, University of Copenhagen, Universitetsparken 5, 2100 Copenhagen \O, Denmark}

\begin{document}

\maketitle

\begin{abstract}
In 2019, Aterias et al. constructed pairs of quantum isomorphic, non-isomorphic graphs from linear constraint systems. This article deals with quantum automorphisms and quantum isomorphisms of colored versions of those graphs. We show that the quantum automorphism group of such a colored graph is the dual of the homogeneous solution group of the underlying linear constraint system. 
Given a vertex-- and edge-colored graph with certain properties, we construct an uncolored graph that has the same quantum automorphism group as the colored graph we started with.
Using those results, we obtain the first known example of a graph that has quantum symmetry and finite quantum automorphism group. Furthermore, we construct a pair of quantum isomorphic, non-isomorphic graphs that both have no quantum symmetry. 
\end{abstract}

\section{Introduction}

Generalizing Mermin's magic square \cite{Mermin}, linear system games were introduced by Cleve and Mittal in \cite{CM}. It was then shown by Cleve, Liu and Slofstra \cite{CLS}, that perfect quantum strategies for those games are related to the representations of the associated solution group; a finitely presented group with generators and relations reflected by the underlying linear system. 

Together with his collaborators, the first author \cite{nonsignalling} constructed quantum isomorphic, but non-isomorphic graphs from linear system games  having perfect quantum strategies but no such classical strategies. In this case, a quantum isomorphism corresponds to a perfect quantum strategy of the isomorphism game which was also introduced in \cite{nonsignalling}. %Magic unitary, quantum automorphisms

Quantum automorphism groups of graphs were first introduced by Banica \cite{QBan} and Bichon \cite{QBic} to obtain examples of quantum permutation groups. They are generalizations of the automorphism group of a graph in the framework of Woronowicz' compact matrix quantum groups. An interesting question to ask is which quantum permutation groups can be realized as the quantum automorphism group of some graph. This question has for example been considered in \cite{QFrucht}. 

In this work, we give a correspondence between representations of solution groups and quantum isomorphisms/quantum automorphisms of colored versions of the graphs appearing in \cite{nonsignalling}. More specifically, we will see that  the quantum automorphism groups of the colored versions of those graphs are given by the dual of the homogeneous solution group of the underlying linear constraint system. Furthermore, representations of the non-homogeneous solution group provide quantum isomorphisms between the colored versions of the graphs. 

Additionally, we will discuss a decoloring procedure for vertex-- and edge-colored graphs which does not change the quantum automorphism group for certain graphs. This procedure decolors the vertices by adding paths of different lengths to vertices of different colors. In a second step, we get rid of the edge-colors by subdividing the colored edges and then adding paths of different lengths to the newly added vertices. 

The previous results allow us to construct two explicit examples of graphs that were not known before: First, we construct a graph that has quantum symmetry and finite quantum automorphism group. Second, we obtain a pair of graphs that are quantum isomorphic and non-isomorphic, where both graphs additionally do not have quantum symmetry. 

%Explain sections 
The article is structured as follows. In Section \ref{secprelim}, we briefly discuss compact quantum groups and give the definition of the quantum automorphism group of a vertex-- and edge-colored graph. In Section \ref{sec:dualsolgroup}, we define colored versions of the graphs appearing in \cite{nonsignalling}. Here we prove that the dual of the solution group is the quantum automorphism group of such a graph. Section \ref{sec:decolor} deals with decoloring the graphs from the previous section without changing the quantum automorphism group. In Section \ref{sec:qsymfinaut}, we provide (uncolored) graphs whose quantum automorphism group is the dual of a solution group. We in particular obtain a graph with quantum symmetry and finite quantum automorphism group, see Corollary \ref{cor:qsymfinqaut}. Finally, we discuss quantum isomorphisms of the colored graphs in Section \ref{sec:qiso}. Here, we present our example of a pair of graphs that are quantum isomorphic but non-isomorphic, where both graphs do not have quantum symmetry, see Corollary \ref{cor:qisonoqsymnocolor}.

\section{Preliminaries}\label{secprelim}

We start with the definition of a compact quantum group, see \cite{CQG}. Throughout this article, we write $A \otimes B$ for the minimal tensor product of the $C^*$-algebras $A$ and $B$.

\begin{definition}
A \emph{compact quantum group} $\mathbb{G}$ is a pair $(C(\mathbb{G}), \Delta)$, where $C(\mathbb{G})$ is a unital $C^*$-algebra and $\Delta: C(\mathbb{G}) \to C(\mathbb{G}) \otimes C(\mathbb{G})$ is a unital $*$-homomorphism with the following properties
\begin{itemize}
\item $(\Delta \otimes \mathrm{id}) \circ \Delta = (\mathrm{id} \otimes \Delta) \circ \Delta$ (coassociativity)
\item $\Delta(C(\mathbb{G}))(1 \otimes C(\mathbb{G}))$ and $\Delta(C(\mathbb{G}))(C(\mathbb{G}) \otimes 1)$ are dense in $C(\mathbb{G}) \otimes C(\mathbb{G})$.
\end{itemize}
\end{definition}

\begin{definition}
Let $\mathbb{G}=(C(\mathbb{G}), \Delta_{\mathbb{G}})$ and $\mathbb{H}=(C(\mathbb{H}), \Delta_{\mathbb{H}})$ be compact quantum groups. Then $\mathbb{G}$ and $\mathbb{H}$ are isomorphic as compact quantum groups if there exists a $*$-isomorphism $\varphi:C(\mathbb{G})\to C(\mathbb{H})$ such that $\Delta_{\mathbb{H}}\circ \varphi=(\varphi \otimes \varphi)\circ \Delta_{\mathbb{G}}$. We then write $\mathbb{G}\cong \mathbb{H}$.
\end{definition}

\begin{definition}
Let $\mathbb{G}=(C(\mathbb{G}), \Delta_{\mathbb{G}})$ be a compact quantum group. We say that $\mathbb{G}$ is \emph{finite} if the $C^*$-algebra $C(\mathbb{G})$ is finite-dimensional. 
\end{definition}

The following example can for example be found in \cite[Example 4.5]{MD}.

\begin{example}\label{Dualqgroup}
Let $\Gamma$ be a discrete group. Note that the group $C^*$-algebra $C^*(\Gamma)$ can be written as 
\begin{align*}
    C^*(\Gamma)=C^*(u_g, g\in \Gamma\,|\, u_g \text{ unitary, }u_{gh}=u_gu_h, u_{g^{-1}}=u^*_g),   
\end{align*}
see for example \cite[A.6]{MW}.
Let $\Delta$ be the $*$-homomorphism $\Delta: C^*(\Gamma) \to C^*(\Gamma) \otimes C^*(\Gamma)$, $\Delta(u_g) = u_g \otimes u_g$. Then the pair $(C^*(\Gamma), \Delta)$ is a compact quantum group known as the \emph{dual} of $\Gamma$ which we denote by $\hat{\Gamma}$.%Note that if \Gamma abelian, then C^*(\Gamma)=C(\hat{\Gamma}), where \hat{\Gamma} is the dual group.
\end{example}

We will now define compact matrix quantum groups, a subclass of compact quantum groups.

\begin{definition}[\cite{CMQG1}]\label{CMQG}
A \emph{compact matrix quantum group} $\mathbb{G}$ is a pair $(C(\mathbb{G}),u)$, where $C(\mathbb{G})$ is a unital $C^*$-algebra and  $u = (u_{ij}) \in M_n(C(\mathbb{G}))$ is a matrix such that 
\begin{itemize}
    \item the elements $u_{ij}$, $1 \le i,j \le n$, generate $C(\mathbb{G})$, 
    \item the $*$-homomorphism $\Delta: C(\mathbb{G}) \to C(\mathbb{G}) \otimes C(\mathbb{G})$, $u_{ij} \mapsto \sum_{k=1}^n u_{ik} \otimes u_{kj}$ exists,
    \item the matrix $u$ and its transpose $u^{T}$ are invertible. 
\end{itemize}
The matrix $u$ is usually called \emph{fundamental representation} of $\mathbb{G}$.\end{definition}

A very important example is the quantum symmetric group, the quantum analogue of the symmetric group. It was defined by Wang in \cite{WanSn}.

\begin{definition}
The \emph{quantum symmetric group} $S_n^+= (C(S_n^+),u)$ is the compact matrix quantum group, where
\begin{align*}
C(S_n^+) := C^*(u_{ij}, \, 1 \le i,j \le n \, | \, u_{ij} = u_{ij}^* = u_{ij}^2, \, \sum_{k=1}^n u_{ik} = \sum_{k=1}^n u_{kj} = 1). 
\end{align*} 
\end{definition}
Note that a matrix $u = (u_{ij})_{i,j \in [n]}$ with entries from a \emph{nontrivial} unital $C^*$-algebra satisfying $u_{ij} = u_{ij}^* = u_{ij}^2$ and $\sum_{k=1}^n u_{ik} = \sum_{k=1}^n u_{kj} = 1$, as in the definition above, is known as a \emph{magic unitary}.

Quantum automorphism groups of finite graphs were defined in \cite{QBan}, \cite{QBic}. We give a more general definition in the following, for vertex-- and edge-colored graphs. Note that edge-colorings are similar to distances in finite quantum metric spaces, as considered in \cite{BanMetric}. For us, graphs are undirected and do not have loops nor multiple edges. A \emph{colored graph} is a graph $G$ with vertex set $V$ and edge set $E$ along with a \emph{coloring function} $c: V \cup E \to S$ for some set $S$. We refer to $c(x)$ as the color of the vertex/edge $x$. To be specific, we will sometimes refer to colored graphs as \emph{vertex-} and \emph{edge-colored} graphs. Furthermore, we will also consider \emph{edge-colored} graphs where the coloring function is defined only on the edge set $E$, i.e., edges but not vertices are colored. In either case we will use $E_c$ to refer to the set of edges of color $c$. Also, $A_{G_c}$ will denote the adjacency matrix of the edge color $c$, i.e., $(A_{G_c})_{ij} = 1$ if $(i,j) \in E_c$ and $(A_{G_c})_{ij} = 0$ otherwise.

\begin{definition}
Let $G$ be a vertex-- and edge-colored graph. The \emph{quantum automorphism group} $\qut(G)$ is the compact matrix quantum group $(C(\qut(G)), u)$, where $C(\qut(G))$ is the universal $C^*$-algebra with generators $u_{ij}$, $i,j \in V(G)$ and relations
\begin{align}
    &u_{ij}=u_{ij}^*=u_{ij}^2, &&i,j \in V(G),\label{rel1}\\
    &\sum_{l} u_{il}= 1 = \sum_{l} u_{li}, &&i\in V(G),\label{rel2}\\
    &u_{ij}=0, &&\text{for all  } i,j \in V(G) \text{ with } c(i)\neq c(j),\label{rel3}\\
    &A_{G_c}u=uA_{G_c} &&\text{for all edge colors $c$,}\label{adjmatrix}
\end{align}
where \eqref{adjmatrix} is nothing but $\sum_k u_{ik}(A_{G_c})_{kj}=\sum_k (A_{G_c})_{ik}u_{kj}$ for all $i,j \in V(G)$ and all edge-colors $c$.
\end{definition}

It is not immediately obvious that this defines a compact matrix quantum group. By Definition \ref{CMQG} we have to show that $u$ and $u^T$ are invertible as well as that the comultiplication $\Delta:C(\qut(G)) \to C(\qut(G))\otimes C(\qut(G)), u_{ij} \mapsto \sum_k u_{ik}\otimes u_{kj}$ is a $*$-homomorphism. It is enough to show $\Delta(\sum_k u_{ik}(A_{G_c})_{kj})=\Delta(\sum_k (A_{G_c})_{ik}u_{kj})$ for $i,j \in V(G)$, $c$ edge-color and $\Delta(u_{ij})=0$ for $c(i)\neq c(j)$, since we know that $S_n^+$ is a compact matrix quantum group. The first equation follows for example from \cite[Lemma 2.1.2]{Schmidtthesis}. For the second equation, we take vertices $i,j$ with $c(i)\neq c(j)$. We have $\Delta(u_{ij})=\sum_k u_{ik} \otimes u_{kj}$. Note that there is no $k \in V(G)$ with $c(k)=c(i)$ and $c(k)=c(j)$, since we assumed $c(i)\neq c(j)$. We deduce $u_{ik}\otimes u_{kj}=0$ for all $k$ and thus $\Delta(u_{ij})=0$. Summarizing, we see that $\qut(G)$ is indeed a compact matrix quantum group.

%\dnote{Maybe we should also point out that relation \ref{rel3} can be written in matrix form.}

The following lemma gives relations that are equivalent to Relation \eqref{adjmatrix}, see \cite[Proposition 2.1.3]{Schmidtthesis}.

\begin{lemma}\label{lem:productrelation}
Let $u_{ij}$, $1 \le i,j \le n,$ be the generators of $C(S_n^+)$. Then Relation \eqref{adjmatrix} is equivalent to the relations
\begin{align*}
    u_{ij}u_{kl}&=0 \text{  if  } (i,k) \in E_c, (j,l)\notin E_c \text{ for some edge-color $c$},\\
    u_{ij}u_{kl}&=0 \text{  if  } (i,k) \notin E_c, (j,l)\in E_c \text{ for some edge-color $c$}.
\end{align*}
\end{lemma}

\section{Colored graphs whose quantum automorphism group is the dual of a solution group}\label{sec:dualsolgroup}

In the following definition we use $1$ to denote the identity element of a group.

\begin{definition}\label{def:solgroup}
Let $M \in \mathbb{F}_2^{m \times n}$ and $b \in \mathbb{F}_2^m$ with $b \ne 0$. The \emph{solution group} $\Gamma(M,b)$ of the linear system $Mx = b$ is the group generated by elements $x_i$ for $i \in [n]$ and an element $\gamma$ satisfying the following relations:
\begin{enumerate}
    \item $x_i^2 = 1$ for all $i \in [n]$;
    \item $x_ix_j = x_jx_i$ if there exists $k \in [m]$ s.t. $M_{ki} = M_{kj} = 1$;
    \item $\prod_{i: M_{ki} = 1} x_i = \gamma^{b_k}$ for all $k \in [m]$;
    \item $\gamma^2 = 1$;
    \item $x_i\gamma = \gamma x_i$ for all $i \in [n]$.
\end{enumerate}

If $b = 0$, then we refer to $\Gamma(M,b)$ as the \emph{homogeneous solution group} of the system $Mx = 0$, and define this the same as above except that we add the relation $\gamma = 1$. This is equivalent to removing $\gamma$ from the list of generators, changing the righthand side of the equation in (3) to $1$, and removing items (4) and (5). We will sometimes also use $\Gamma_0(M)$ to denote this group.

We will denote by $S_k(M)$ the set $\{i \in [n] : M_{ki} = 1\}$, often writing simply $S_k$ when $M$ is clear from context. We use $\pm 1^S$ to denote the set of functions $\alpha : S \to \{1,-1\}$, and will typically write $\alpha_i$ instead of $\alpha(i)$. We will also use $\pm 1^S_0$ to denote the subset of such functions satisfying $\prod_{i \in S} \alpha_i = 1$, and similarly use $\pm 1^S_1$ for the set of such functions satisfying $\prod_{i \in S} \alpha_i = -1$.
\end{definition}

\begin{definition}\label{defgraph}
Let $M \in \mathbb{F}_2^{m \times n}$ and $b \in \mathbb{F}_2^m$. Define the colored graph $G := G(M,b)$ as follows. The vertex set of $G$ is $\left\{(k,\alpha): k \in [m], \ \alpha \in \pm 1^{S_k}_{b_k}\right\}$. The color of a vertex $v = (k,\alpha)$, denoted $c(v)$, is $k$. We also define the subsets $V_k = \{v \in V(G) : c(v) = k\}$ for all $k \in [m]$, and these partition the vertex set. For any $l,k \in [m]$ such that $S_l \cap S_k \ne \varnothing$, the graph $G$ contains all (non-loop) edges between $V_l$ and $V_k$ (thus each $V_k$ induces a complete subgraph). Given an edge $e$ between adjacent vertices $(l,\alpha)$ and $(k,\beta)$, the color of $e$, denoted $c(e)$, is equal to the function $\alpha \triangle \beta \in \pm 1^{S_l \cap S_k}$ defined as $(\alpha \triangle \beta)_i = \alpha_i\beta_i$.
\end{definition}

Note that $\alpha\triangle \beta = \beta \triangle \alpha$ and so the edge colors defined above really are edge colors and not arc (directed edge) colors.

\begin{remark}\label{rem1}
According to the above definition, it is possible for two edges $e$ and $e'$ between pairs of vertices $(l,\alpha),(k,\beta)$ and $(l',\alpha'),(k',\beta')$ to be colored the same color even if $l$ is not equal to either $l'$ nor $k'$, i.e., the edges are between different pairs of subsets $V_l,V_k$ and $V_{l'},V_{k'}$. This can happen since it is possible that $S_l \cap S_k = S_{l'} \cap S_{k'}$ and $\alpha \triangle \beta = \alpha' \triangle \beta'$. However, we do wish such edges to be distinguished by the (quantum) automorphism group of $G(M,b)$, i.e., $u_{(l,\alpha),(l'\alpha')}u_{(k,\beta),(k'\beta')} = 0$ where $u$ is the fundamental representation of $\qut(G(M,b))$\footnote{This is the quantum analog of there being no automorphism mapping the edge $(l,\alpha)(k,\beta)$ to $(l'\alpha')(k',\beta')$.}. This could be done explicitly by defining the color of the edge $e$ (for instance) to be $c(e) = (\{l,k\},\alpha \triangle \beta)$. However, this is redundant for our purposes since such edges $e$ and $e'$ as described above are already distinguished by the (quantum) automorphism group due to the vertex colors of the endpoints of the edges. In other words, the edges between $V_l$ and $V_k$ can be thought of as being \emph{implicitly} colored distinctly from those between $V_{l'}$ and $V_{k'}$ whenever $\{l,k\} \ne \{l',k'\}$.
\end{remark}

\begin{remark}\label{rem2}
We also note that in order to reduce the total number of colors used, for each pair $l,k \in [m]$, we can choose one color $\alpha \in \pm 1^{S_l \cap S_k}$ and replace these edges with non-edges. Moreover, instead of coloring the edges between $V_l$ and $V_k$ with functions from $\pm 1^{S_l \cap S_k}$, we can simply use $\{1, \ldots, 2^{|S_l\cap S_k|}-1\}$, and this will not change $\qut(G)$ as explained in the previous remark. In fact we will need to do this for some of our results later on.
\end{remark}

\begin{example}
Consider the linear system $Mx=b$, where \begin{align*}
 M=\begin{pmatrix}1&1&1&0&0\\1&0&0&1&1\end{pmatrix},b=\begin{pmatrix}0\\1\end{pmatrix}.   
\end{align*}
Then the graph $G(M,b)$ is the one given in Figure \ref{figure1}. 
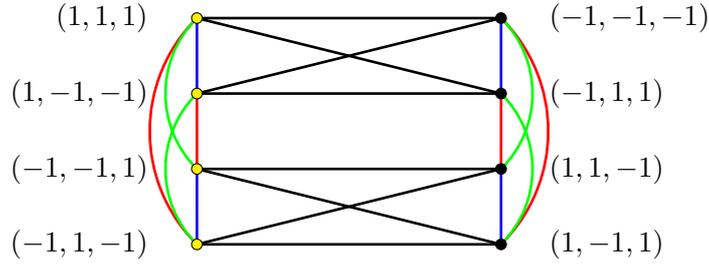
\begin{figure}[h]
\begin{center}\begin{tikzpicture}[scale=1, baseline=-1.6cm]
%\coordinate[label=left:$G\text{ =}$] (1) at (180:1cm);
\draw[blue, line width=1pt] (-2,-1.5) -- (-2,-0.5);
\draw[blue, line width=1pt] (-2,1.5) -- (-2,0.5);
\draw[blue, line width=1pt] (2,-1.5) -- (2,-0.5);
\draw[blue, line width=1pt] (2,1.5) -- (2,0.5);
\draw[red, line width=1pt] (2,-0.5) -- (2,0.5);
\draw[red, line width=1pt] (-2,-0.5) -- (-2,0.5);
\draw[red, line width=1pt] (-2,-1.5) to[in=-135,out=-225] (-2,1.5);
\draw[red, line width=1pt] (2,-1.5) to[in=-45,out=45] (2,1.5);
\draw[green, line width=1pt] (2,-1.5) to[in=-45,out=45] (2,0.5);
\draw[green, line width=1pt] (2,-0.5) to[in=-45,out=45] (2,1.5);
\draw[green, line width=1pt] (-2,-1.5) to[in=-135,out=-225] (-2,0.5);
\draw[green, line width=1pt] (-2,-0.5) to[in=-135,out=-225] (-2,1.5);
\draw[line width=1pt] (-2,1.5) -- (2,1.5);
\draw[line width=1pt] (-2,1.5) -- (2,0.5);
\draw[line width=1pt] (-2,0.5) -- (2,1.5);
\draw[line width=1pt] (-2,0.5) -- (2,0.5);
\draw[line width=1pt] (-2,-0.5) -- (2,-1.5);
\draw[line width=1pt] (-2,-0.5) -- (2,-0.5);
\draw[line width=1pt] (-2,-1.5) -- (2,-1.5);
\draw[line width=1pt] (-2,-1.5) -- (2,-0.5);
\foreach \x in {-1.5,-0.5,0.5,1.5}
{\draw[black,fill=black] (2,\x) circle (2pt);}
\foreach \x in {-1.5,-0.5,0.5,1.5}
{\draw[black,fill=yellow] (-2,\x) circle (2pt);}
\coordinate[label=left:{$(1,1,1)$}] (1) at (-2.5,1.5);
\coordinate[label=left:{$(1,-1,-1)$}] (1) at (-2.5,0.5);
\coordinate[label=left:{$(-1,-1,1)$}] (1) at (-2.5,-0.5);
\coordinate[label=left:{$(-1,1,-1)$}] (1) at (-2.5,-1.5);
\coordinate[label=right:{$(-1,-1,-1)$}] (1) at (2.5,1.5);
\coordinate[label=right:{$(-1,1,1)$}] (1) at (2.5,0.5);
\coordinate[label=right:{$(1,1,-1)$}] (1) at (2.5,-0.5);
\coordinate[label=right:{$(1,-1,1)$}] (1) at (2.5,-1.5);
%\coordinate[label=left:{$x_1x_2x_3=1$}] (1) at (-5,0);
%\coordinate[label=right:{$x_1x_4x_5=-1$}] (1) at (5,0);
\end{tikzpicture}\end{center}\caption{The graph $G(M,b)$ for $M$ and $b$ as above. The vertices on the left hand side are the solutions of $x_1x_2x_3=1$, the vertices on the right hand side are the solutions of $x_1x_4x_5=-1$. We used Remarks \ref{rem1} and \ref{rem2} to reduce the number of colors needed in the graph.}\label{figure1}
\end{figure}\end{example}

\begin{lemma}\label{lem:qutGdecomp}
Let $M \in \mathbb{F}_2^{m \times n}$ and $b \in \mathbb{F}_2^m$. Set $G = G(M,b)$. If $u$ is the fundamental representation of $\qut(G)$, then $u$ has the form
%\[u = \bigoplus_{k \in [m]} u^{(k)},\]
\begin{align}
\begin{blockarray}{cccccc}
&V_1&V_2&\dots&\dots&V_m\\
\begin{block}{c(ccccc)}
V_1&u^{(1)}&0&0&\dots&0\\ V_2&0&u^{(2)}&0&\dots&0\\ \vdots&0&0&u^{(3)}&\dots&0\\\vdots&\vdots &\vdots&\vdots&\ddots&0\\V_m&0&0&0&0&u^{(m)} \\
\end{block}
\end{blockarray}
\end{align}
%where each $u^{(k)}$ is a magic unitary indexed by the subset $V_k \subseteq V(G)$. 
where each $u^{(k)}$ is a magic unitary. Furthermore, $u^{(k)}_{\alpha,\beta} := u_{(k,\alpha),(k,\beta)}$ depends only on $k$ and the value of $\alpha \triangle \beta$ and therefore the entries of $u^{(k)}$ pairwise commute for each $k \in [m]$.
\end{lemma}
\proof
If $l \ne k$, then $c((l,\alpha)) = l \ne k = c((k,\beta))$ and thus $u_{(l,\alpha),(k,\beta)} = 0$. This shows that $u$ has the block form given in the lemma statement. Moreover, each diagonal block must be a magic unitary since $u$ is a magic unitary.

Now fix $k \in [m]$. Pick $\tilde{\alpha},\tilde{\beta} \in \pm 1^{S_k}_{b_k}$ arbitrarily. Suppose that $\alpha,\beta \in \pm 1^{S_k}_{b_k}$ are such that $\alpha \triangle \beta = \tilde{\alpha} \triangle \tilde{\beta}$. Note that since all these functions are elements of $\pm 1^{S_k}_{b_k}$, the operation $\triangle$ is simply the pointwise product, and therefore we have that $\tilde{\alpha} \triangle \alpha = \tilde{\beta} \triangle \beta$. In particular, this implies that the color of the edge between $(k,\tilde{\alpha})$ and $(k,\alpha)$ is equal to the color of the edge between $(k,\tilde{\beta})$ and $(k,\beta')$ if and only if $\beta' = \beta$. Therefore, $u^{(k)}_{\tilde{\alpha},\tilde{\beta}}u^{(k)}_{\alpha,\beta'} = 0$ unless $\beta' = \beta$. We similarly have that $u^{(k)}_{\tilde{\alpha},\tilde{\beta}'}u^{(k)}_{\alpha,\beta} = 0$ unless $\tilde{\beta}' = \tilde{\beta}$. It follows that
\[u^{(k)}_{\tilde{\alpha},\tilde{\beta}} = u^{(k)}_{\tilde{\alpha},\tilde{\beta}}\left(\sum_{\beta' \in \pm 1^{S_k}_{b_k}} u^{(k)}_{\alpha,\beta'}\right) = u^{(k)}_{\tilde{\alpha},\tilde{\beta}} u^{(k)}_{\alpha,\beta} = \left(\sum_{\tilde{\beta}' \in \pm 1^{S_k}_{b_k}} u^{(k)}_{\tilde{\alpha},\tilde{\beta}'}\right) u^{(k)}_{\alpha,\beta} = u^{(k)}_{\alpha,\beta}.\]
Thus the operator $u^{(k)}_{\alpha,\beta}$ depends only on the value of $\alpha \triangle \beta$ (and $k$) as claimed. Let us therefore use $u^{(k)}_\delta$ to denote $u^{(k)}_{\alpha,\beta}$ where $\delta = \alpha \triangle \beta$. For any fixed $\alpha \in \pm 1^{S_k}_{b_k}$, the set $\{\alpha \triangle \beta : \beta \in \pm 1^{S_k}_{b_k}\}$ is equal to $\pm 1^{S_k}_0$. Thus the set of operators appearing in any row (and similarly any column) of $u^{(k)}$ is $\{u^{(k)}_\delta : \delta \in \pm 1^{S_k}_{0}\}$. Since the entries of any given row of a magic unitary commute, and every row of $u^{(k)}$ contains the same set of entries, we have that all entries of $u^{(k)}$ commute.\qeds

%\dnote{In the following we use the fact that $C^*(\Gamma)$ for a group $\Gamma$ can be expressed as a universal $C^*$-algebra in a particular way. So we should include an explanation of this somewhere.}
\begin{remark}\label{rem:unigrpalg}
Let $\Gamma=\Gamma(M,0)$ be the homogeneous solution group of the system $Mx=0$. By Example \ref{Dualqgroup}, the group $C^*$-algebra $C^*(\Gamma)$ can be written as 
\begin{align*}
    C^*(\Gamma)=C^*(x_i, i \in [n]\,|\, &x_i=x_i^*, x_i^2 = 1, x_ix_j = x_jx_i \text{ if there exists $k \in [m]$ s.t. } M_{ki} = M_{kj} = 1,\\&\prod_{i: M_{ki} = 1} x_i = 1\text{ for all $k \in [m]$}).
\end{align*}
Note that we slightly abuse the notation by writing $x_i$ instead of $u_{x_i}$.
\end{remark}

\begin{theorem}\label{thm:qaut2cgamma}
Let $M \in \mathbb{F}_2^{m \times n}$ and $b \in \mathbb{F}_2^m$. Set $G = G(M,b)$ and $\Gamma = \Gamma(M,0) = \Gamma_0(M)$. Let $\Delta_\Gamma : C^*(\Gamma) \to C^*(\Gamma) \otimes C^*(\Gamma)$ denote the coproduct given by $\Delta_\Gamma(g) = g \otimes g$ for $g \in \Gamma$, and let $\Delta_G: C(\qut(G)) \to C(\qut(G)) \otimes C(\qut(G))$ denote the coproduct given by $\Delta_G(u_{a,c}) = \sum_{b \in V(G)} u_{a,b} \otimes u_{b,c}$ where $u = (u_{a,c})$ is the fundamental representation of $\qut(G)$. Then there exists an isomorphism $\varphi : C^*(\Gamma) \to C(\qut(G))$ such that $\Delta_G \circ \varphi = (\varphi \otimes \varphi) \circ \Delta_\Gamma$.
\end{theorem}
\proof
The proof is structured as follows: we use the universal properties of $C^*(\Gamma)$ and $C(\qut(G))$ respectively to first obtain a $*$-homomorphism $\varphi_1: C^*(\Gamma) \to C(\qut(G))$, and then obtain a $*$-homomorphism $\varphi_2: C(\qut(G)) \to C^*(\Gamma)$. We will then show that $\varphi_1$ and $\varphi_2$ are inverses of each other thus proving that they are in fact isomorphisms. Lastly, we will prove that $\varphi_1$ intertwines the coproducts as described in the theorem statement.\\

\noindent\emph{Step 1: Construction of a $*$-homomorphism $\varphi_1: C^*(\Gamma) \to C(\qut(G))$.}\\
For this step, we will construct elements $y_i$ of $C(\qut(G))$ that satisfy the relations of the generators $x_i$ of $C^*(\Gamma)$ given in Remark~\ref{rem:unigrpalg}. This proves that there is a $*$-homomorphism $\varphi_1$ from $C^*(\Gamma)$ to $C(\qut(G))$ such that $\varphi_1(x_i) = y_i$. Later we will see that $\varphi_1$ is in fact an isomorphism.

Let $u$ be the fundamental representation of $\qut(G)$. By Lemma~\ref{lem:qutGdecomp}, $u = \bigoplus_{k \in [m]} u^{(k)}$ where each $u^{(k)}$ is a magic unitary indexed by $V_k = \{(k,\alpha) : \alpha \in \pm 1^{S_k}_{b_k}\}$. Moreover, the operator $u^{(k)}_{\alpha,\beta} := u_{(k,\alpha),(k,\beta)}$ depends only on $k$ and the value of $\alpha \triangle \beta \in \pm 1^{S_k}_{0}$. Thus, as in the proof of Lemma~\ref{lem:qutGdecomp}, for each $\delta \in \pm 1^{S_k}_0$ we let $u^{(k)}_\delta$ denote $u^{(k)}_{\alpha,\beta}$ such that $\alpha \triangle \beta = \delta$. Note that for any $\alpha \in \pm 1^{S_k}_{b_k}$, we have that the set $\{\alpha \triangle \beta : \beta \in \pm 1^{S_k}_{b_k}\}$ is equal to $\pm 1^{S_k}_0$ and thus every row/column of $u^{(k)}$ contains the same set of operators and
\begin{equation}\label{eq:udeltares}
    \sum_{\delta \in \pm 1^{S_k}_0} u^{(k)}_\delta = \sum_{\beta \in \pm 1^{S_k}_{b_k}} = \sum_{\alpha \in \pm 1^{S_k}_{b_k}} u^{(k)}_{\alpha,\beta} = 1.
\end{equation}

Now we define $y^{(k)}_i = \sum_{\alpha \in \pm 1^{S_k}_0} \alpha_i u^{(k)}_\alpha$ for all $k \in [m]$ and $i \in S_k$. We first aim to show that $y^{(k)}_i$ does not depend on $k$.

Fix $i \in [n]$ and $l,k \in [m]$ such that $i \in S_l \cap S_k$. Consider the edges between the subsets $V_l$ and $V_k$. For each $\delta \in \pm 1^{S_l \cap S_k}$, let $A^\delta$ be the adjacency matrix of the graph consisting of the edges of $G$ colored $\delta$. Further, let $B^\delta$ be the submatrix of $A^\delta$ consisting of the rows indexed by $V_l$ and columns indexed by $V_k$. In other words,
\[B^\delta_{(l,\alpha),(k,\beta)} = \begin{cases}1 & \text{if } \alpha \triangle \beta = \delta \\ 0 & \text{o.w.}\end{cases}\]
Also define $A^+ = \sum_{\delta \in \pm 1^{S_l \cap S_k}, \delta_i = +1} A^\delta$ and $A^- = \sum_{\delta \in \pm 1^{S_l \cap S_k}, \delta_i = -1} A^\delta$, and similarly define $B^+$ and $B^-$. Then
\[B^+_{(l,\alpha),(k,\beta)} = \begin{cases}1 & \text{if } \alpha_i = \beta_i \\ 0 & \text{o.w.}\end{cases} \quad \text{ \& } \quad B^-_{(l,\alpha),(k,\beta)} = \begin{cases}1 & \text{if } \alpha_i = -\beta_i \\ 0 & \text{o.w.}\end{cases}\]
Since $u$ must commute with each $A^\delta$ by definition of $\qut(G)$, it must also commute with both $A^+$ and $A^-$. This implies that $u^{(l)}B^+ = B^+ u^{(k)}$ and $u^{(l)}B^- = B^- u^{(k)}$. Considering the $(l,\alpha),(k,\beta)$ entry of both sides of the former equation, we see that
\[\sum_{\substack{\alpha' \in \pm 1^{S_l}_{b_l} \\ \text{s.t. } \alpha'_i =  \beta_i}}u^{(l)}_{\alpha,\alpha'} = \sum_{\substack{\beta' \in \pm 1^{S_k}_{b_k} \\ \text{s.t. } \beta'_i = \alpha_i}}u^{(k)}_{\beta',\beta}.\]
Note that for every term in the above sums, we have that $\alpha_i\alpha'_i = \alpha_i\beta_i = \beta'_i\beta_i$. Thus, reexpressing $u^{(l)}_{\alpha,\alpha'}$ as $u^{(l)}_{\tilde{\alpha}}$ where $\tilde{\alpha} = \alpha \triangle \alpha'$, and similarly for $u^{(k)}_{\beta',\beta}$ we have that
\[\sum_{\substack{\tilde{\alpha} \in \pm 1^{S_l}_{0} \\ \text{s.t. } \tilde{\alpha}_i =  \alpha_i\beta_i}}u^{(l)}_{\tilde{\alpha}} = \sum_{\substack{\tilde{\beta} \in \pm 1^{S_k}_{0} \\ \text{s.t. } \tilde{\beta}_i = \alpha_i\beta_i}}u^{(k)}_{\tilde{\beta}}.\]
Doing the same for $u^{(l)}B^- = B^- u^{(k)}$ yields the same equation but with $\alpha_i\beta_i$ replaced with $-\alpha_i\beta_i$ and combining these proves that
\[\sum_{\tilde{\alpha} \in \pm 1^{S_l}_{0}, \tilde{\alpha}_i = x}u^{(l)}_{\tilde{\alpha}} = \sum_{\tilde{\beta} \in \pm 1^{S_k}_{0}, \tilde{\beta}_i = x}u^{(k)}_{\tilde{\beta}}\]
for any $x \in \{+1,-1\}$. Therefore,
\[y^{(l)}_i = \left(\sum_{\substack{\alpha \in \pm 1^{S_l}_0 \\ \text{s.t. } \alpha_i = +1}} u^{(l)}_{\alpha}\right) - \left(\sum_{\substack{\alpha \in \pm 1^{S_l}_0 \\ \text{s.t. } \alpha_i = -1}} u^{(l)}_{\alpha}\right) = \left(\sum_{\substack{\alpha \in \pm 1^{S_k}_0 \\ \text{s.t. } \alpha_i = +1}} u^{(k)}_{\alpha}\right) - \left(\sum_{\substack{\alpha \in \pm 1^{S_k}_0 \\ \text{s.t. } \alpha_i = -1}} u^{(k)}_{\alpha}\right) = y^{(k)}_i.\]
So we have shown that the value of $y^{(k)}_i$ does not depend on $k$, and thus we will simply denote this operator by $y_i$.

Now note that since $y_i$ is a linear combination (with real coefficients) of the operators $u^{(k)}_\alpha$ for $k \in [m]$ such that $i \in S_k$, and these operators are entries of the magic unitary $u$, we have that $y_i^* = y_i$. Also, by Equation~\eqref{eq:udeltares} we have that $u^{(k)}_\alpha u^{(k)}_\beta = 0$ for $\alpha \ne \beta$ and thus
\[y_i^2 = \left(\sum_{\alpha \in \pm 1^{S_k}_0} \alpha_i u^{(k)}_\alpha\right)^2 = \sum_{\alpha \in \pm 1^{S_k}_0} \alpha^2_i \left(u^{(k)}_\alpha\right)^2 = \sum_{\alpha \in \pm 1^{S_k}_0} u^{(k)}_\alpha = 1.\]
Thus the $y_i$ satisfy relation (1) from Definition~\ref{def:solgroup}.

Next we will show that relation (2) of Definition~\ref{def:solgroup} holds, i.e., that $y_iy_j = y_jy_i$ if there exists $k \in [m]$ such that $i,j \in S_k$. Suppose that $i,j,k$ are as described. Then $y_i = y_i^{(k)}$ and $y_j = y_j^{(k)}$ are both linear combinations of the entries of $u^{(k)}$ which pairwise commute by Lemma~\ref{lem:qutGdecomp}. Therefore $y_i$ and $y_j$ commute as desired.

Lastly, we must show that relation (3) of Definition~\ref{def:solgroup} holds. Recall that we are trying to show that $C^*(\Gamma_0(M)) \cong C(\qut(G))$, i.e., we are in the \emph{homogeneous} solution group case. Thus we must show that $\prod_{i \in S_k}y_i = 1$ for all $k \in [m]$. We have that
\begin{align*}
    \prod_{i \in S_k}y_i &= \prod_{i \in S_k} y^{(k)}_i \\
    &= \prod_{i \in S_k} \left(\sum_{\alpha \in \pm 1^{S_k}_0} \alpha_i u^{(k)}_\alpha\right)
\end{align*}
Since the elements $u^{(k)}_\alpha$ are pairwise orthogonal, when we expand the above product all cross terms disappear and we obtain
\[\prod_{i \in S_k}y_i = \sum_{\alpha \in \pm 1^{S_k}_0} \left(\prod_{i \in S_k} \alpha_i\right)u^{(k)}_\alpha = \sum_{\alpha \in \pm 1^{S_k}_0} u^{(k)}_\alpha = 1,\]
as desired. Therefore the elements $y_i \in C(\qut(G))$ for $i \in [n]$ satisfy the relations of the generators of $C^*(\Gamma)$ and thus there exists a $*$-homomorphism $\varphi_1$ from $C^*(\Gamma)$ to $C(\qut(G))$ such that $\varphi_1(x_i) = y_i$ for all $i \in [n]$.\\

\noindent\emph{Step 2: Construction of a $*$-homomorphism $\varphi_2: C(\qut(G)) \to C^*(\Gamma)$.}\\
For this step, we will construct elements $v_{(l,\alpha),(k,\beta)}$ of $C^*(\Gamma)$ that satisfy all the relations of the generators $u_{(l,\alpha),(k,\beta)}$ of $C(\qut(G))$. This proves that there is a $*$-homomorphism $\varphi_2$ from $C(\qut(G))$ to $C^*(\Gamma)$ such that $\varphi_2(u_{(l,\alpha),(k,\beta)}) = v_{(l,\alpha),(k,\beta)}$. Later we will see that $\varphi_2$ is an isomorphism.

For each $i \in [n]$, define $p^\pm_i = \frac{1}{2}(1 \pm x_i)$. It is straightforward to check that $p_i^\pm$ is a projection, i.e., $(p^\pm_i)^* = p^\pm_i = (p^\pm_i)^2$, and that $x_i = p_i^+ - p_i^-$ and $1 = p_i^+ + p_i^-$. Note that the latter implies that $p_i^+$ and $p_i^-$ are orthogonal, i.e., $p_i^+p_i^- = 0$. We will abuse notation somewhat and write $p_i^{\alpha_i}$ to denote $p_i^+$ whenever $\alpha_i = +1$ and similarly for $p_i^-$ when $\alpha_i = -1$. Next, for all $k \in [m]$ and $\alpha \in \pm 1^{S_k}$ define $v^{(k)}_\alpha = \prod_{i \in S_k} p_i^{\alpha_i}$. Note that the $v^{(k)}_\alpha$ are indeed well defined since the commutativity of the elements $p^\pm_i$ for $i \in S_k$ follows from the commutativity of the $x_i$ for $i \in S_k$ (i.e., relation (2) from Definition~\ref{def:solgroup}). Since it is the product of pairwise commuting projections, we have that $v^{(k)}_\alpha$ is a projection. Furthermore,
\begin{equation}\label{eq:vres}
    \sum_{\alpha \in \pm 1^{S_k}} v^{(k)}_\alpha = \sum_{\alpha \in \pm 1^{S_k}} \prod_{i \in S_k} p_i^{\alpha_i} = \prod_{i \in S_k} (p^+_i + p^-_i) = \prod_{i \in S_k} 1 = 1.
\end{equation}
This also implies that for a fixed $k \in [m]$ the projections $v^{(k)}_\alpha$ are pairwise orthogonal.

Now suppose that $\alpha \in \pm 1^{S_k}_1$, i.e., that $\prod_{i \in S_k} \alpha_i = -1$. Then, using the easily checked fact that $x_i p_i^{\alpha_i} = \alpha_i p_i^{\alpha_i}$, we have that
\begin{equation}\label{eq:vka0}
v^{(k)}_\alpha = \left(\prod_{i \in S_k} x_i\right) v^{(k)}_\alpha = \prod_{i \in S_k} x_ip_i^{\alpha_i} = \prod_{i \in S_k} \alpha_i p_i^{\alpha_i} = \left(\prod_{i \in S_k} \alpha_i\right) v^{(k)}_\alpha = -v^{(k)}_\alpha.
\end{equation}
Thus $v^{(k)}_\alpha = 0$ for all $\alpha \in \pm 1^{S_k}_1$. Combining this with Equation~\eqref{eq:vres}, we have that
\begin{equation}\label{eq:vres0}
    \sum_{\alpha \in \pm 1^{S_k}_0} v^{(k)}_\alpha = 1.
\end{equation}

Next, for all $k \in [m]$ and $\alpha,\beta \in \pm 1^{S_k}_{b_k}$ we define $v^{(k)}_{\alpha,\beta} := v^{(k)}_\delta$ where $\delta = \alpha \triangle \beta$. Lastly, define $v$ to the the matrix indexed by $V(G)$ such that
\[v_{(l,\alpha),(k,\beta)} = \begin{cases} v^{(l)}_{\alpha, \beta} & \text{if } l = k \\ 0 & \text{o.w.}\end{cases}\]
We aim to show that $v$ satisfies all of the relations satisfied by the fundamental representation $u$ of $\qut(G)$. First, $v$ is a magic unitary since all of its entries are projections and the sum of its $(l,\alpha)$-row is simply $\sum_{\beta \in \pm 1^{S_l}_{b_l}} v^{(l)}_{\alpha,\beta} = \sum_{\delta \in \pm 1^{S_l}_0} v^{(l)}_{\delta} = 1$. Next we must show that for $a,b \in V(G)$, we have $v_{a,b} = 0$ whenever $c(a) \ne c(b)$, i.e., the colors of $a$ and $b$ are different. Recall from the definition of $G(M,b)$ that the color of the vertex $(l,\alpha)$ is $l$. Thus $v$ satisfies this relation by definition. Lastly, we must show that for $a,b,a',b' \in V(G)$, we have $v_{a,b}v_{a',b'} = 0$ whenever the edges $\{a,a'\}$ and $\{b,b'\}$ have different colors (or when one is an edge and the other is not). The previously shown relation already implies this one unless $c(a) = c(b)$ and $c(a') = c(b')$. So we may assume that $a = (l,\alpha)$, $b = (l,\beta)$, $a' = (k,\alpha')$, and $b' = (k,\beta')$. Thus, letting $\delta = \alpha \triangle \beta$ and $\delta' = \alpha' \triangle \beta'$, we have $v_{a,b} = v^{(l)}_{\alpha,\beta} = v^{(l)}_\delta$ and $v_{a',b'} = v^{(k)}_{\alpha',\beta'} = v^{(k)}_{\delta'}$. Suppose the edges $\{a,a'\}$ and $\{b,b'\}$ do have different colors. By definition of $G(M,b)$ this is equivalent to $\alpha \triangle \alpha' \ne \beta \triangle \beta'$. This means there exists $j \in S_l \cap S_k$ such that $\alpha_j\alpha'_j \ne \beta_j\beta'_j$. Since all these terms are in $\{+1,-1\}$, we have that $\delta_j = \alpha_j\beta_j \ne \alpha'_j\beta'_j = \delta'_j$. This implies that $p_j^{\delta_j}p_j^{\delta'_j} = 0$.  Therefore,
\[v_{a,b}v_{a',b'} = v^{(l)}_{\delta}v^{(k)}_{\delta'} = \left(\prod_{i \in S_l} p_i^{\delta_i}\right)\left(\prod_{i \in S_k} p_i^{\delta'_i}\right) = p_j^{\delta_j}p_j^{\delta'_j}\left(\prod_{i \in S_l \setminus \{j\}} p_i^{\delta_i}\right)\left(\prod_{i \in S_k \setminus \{j\}} p_i^{\delta'_i}\right) = 0.\]
So we have shown that $v$ is a magic unitary satisfying all of the relations satisfied by $u$ and thus by the universal property of $C(\qut(G))$, there exists a $*$-homomorphism $\varphi_2$ from $C(\qut(G))$ to $C^*(\Gamma)$ such that $\varphi_2(u_{(l,\alpha),(k,\beta)}) = v_{(l,\alpha),(k,\beta)}$, which is equivalent to $\varphi_2\left(u^{(k)}_\delta\right) = v^{(k)}_\delta$ for all $k \in [m]$ and $\delta \in \pm 1^{S_k}_0$.\\

\noindent\emph{Step 3: Showing that $\varphi_1$ and $\varphi_2$ are inverses of each other.}\\
We first show that $\varphi_2 \circ \varphi_1: C^*(\Gamma) \to C^*(\Gamma)$ is the identity. Since both $\varphi_1$ and $\varphi_2$ are $*$-homomorphisms, it suffices to show that $\varphi_2 \circ \varphi_1$ acts as the identity on the generators $x_i$ of $C^*(\Gamma)$. From Step 1 we have that $\varphi_1(x_i) = y_i = \sum_{\alpha \in \pm 1^{S_k}_0} \alpha_i u^{(k)}_\alpha$ for any $k \in [m]$ such that $i \in S_k$. From Step 2 we have that $\varphi_2(u^{(k)}_\alpha) = v^{(k)}_\alpha = \prod_{i \in S_k} p_i^{\alpha_i}$. Therefore,
\[\varphi_2 \circ \varphi_1 (x_i) = \sum_{\alpha \in \pm 1^{S_k}_0} \alpha_i \varphi_2(u^{(k)}_\alpha) = \sum_{\alpha \in \pm 1^{S_k}_0} \alpha_i v^{(k)}_\alpha.\]
To see that the final expression above is equal to $x_i$, we will show that multiplying it by $x_i$ yields 1. First, note that since $x_ip_i^{\alpha_i} = \alpha_ip_i^{\alpha_i}$, we have that
\[x_iv^{(k)}_\alpha = x_i \prod_{j \in S_k} p_j^{\alpha_j} = x_ip_i^{\alpha_i}\prod_{j \in S_k\setminus \{i\}} p_j^{\alpha_j} = \alpha_ip_i^{\alpha_i}\prod_{j \in S_k\setminus \{i\}} p_j^{\alpha_j} = \alpha_i\prod_{j \in S_k} p_j^{\alpha_j} = \alpha_i v^{(k)}_\alpha.\]
Thus, making use of Equation~\eqref{eq:vres0} in the last equality, we have that
\[x_i\left(\varphi_2 \circ \varphi_1 (x_i)\right) = x_i\left(\sum_{\alpha \in \pm 1^{S_k}_0} \alpha_i v^{(k)}_\alpha\right) = \sum_{\alpha \in \pm 1^{S_k}_0} \alpha_i x_iv^{(k)}_\alpha = \sum_{\alpha \in \pm 1^{S_k}_0} \alpha_i^2 v^{(k)}_\alpha = \sum_{\alpha \in \pm 1^{S_k}_0} v^{(k)}_\alpha = 1,\]
as desired. It follows that $\varphi_2 \circ \varphi_1(x_i)$ is the inverse of $x_i$, which is its own inverse by definition. Therefore $\varphi_2 \circ \varphi_1(x_i) = x_i$ and thus $\varphi_2 \circ \varphi_1$ is the identity map on $C^*(\Gamma)$.

Now we will show that $\varphi_1 \circ \varphi_2 : C(\qut(G)) \to C(\qut(G))$ is the identity map. As above, it suffices to show that $\varphi_1 \circ \varphi_2$ acts as the identity on all entries of the fundamental representation $u$. Other than $0$, the entries of $u$ are precisely the elements $u^{(k)}_\alpha$ for $k \in [m]$ and $\alpha \in \pm 1^{S_k}_0$. From Step 2 we have that $\varphi_2(u^{(k)}_\alpha) = v^{(k)}_\alpha = \prod_{i \in S_k} p_i^{\alpha_i} = \prod_{i \in S_k} \frac{1}{2}(1 + \alpha_i x_i)$, and from Step 1 we have that $\varphi_1(x_i) = y_i = \sum_{\alpha \in \pm 1^{S_k}_0} \alpha_i u^{(k)}_\alpha$. Therefore,
\begin{equation}\label{eq:phiu}
    \varphi_1 \circ \varphi_2(u^{(k)}_\alpha) = \varphi_1 \left(\prod_{i\in S_k} \frac{1}{2}(1 + \alpha_i x_i)\right) = \prod_{i\in S_k} \frac{1}{2}(1 + \alpha_i y_i),
\end{equation}
where we have used the fact that $\varphi_1$ maps $1 \in C^*(\Gamma)$ to $1 \in C(\qut(G))$ which follows from $\varphi_1(1) = \varphi_1(x_i^2) = \varphi_1(x_i)^2 = y_i^2 = 1$. Let us use $\tilde{u}^{(k)}_\alpha$ to denote the final expression of Equation~\eqref{eq:phiu}. We wish to show that $\tilde{u}^{(k)}_\alpha = u^{(k)}_\alpha$. Since the $y_i$ satisfy all of the relations that the $x_i$ satisfy, the same argument as for the $v^{(k)}_\alpha$ can be used to show that, for fixed $k \in [m]$ the $\tilde{u}^{(k)}_\alpha$ are pairwise orthogonal projections satisfying $\sum_{\alpha \in \pm 1^{S_k}_0} \tilde{u}^{(k)}_\alpha = 1$. For $i \in S_k$ and $\beta \in \pm 1^{S_k}_0$ we have that
\[u^{(k)}_\beta y_i = u^{(k)}_\beta\sum_{\alpha \in \pm 1^{S_k}_0} \alpha_i u^{(k)}_\alpha = \beta_i u^{(k)}_\beta.\]
Suppose that $\alpha \ne \beta$. Then there exists $j \in S_k$ such that $\alpha_j \ne \beta_j$ and therefore
\begin{align*}
    u^{(k)}_\beta \tilde{u}^{(k)}_\alpha &= u^{(k)}_\beta \prod_{i \in S_k} \frac{1}{2}(1 + \alpha_iy_i) \\
    &= u^{(k)}_\beta \frac{1}{2}(1 + \alpha_j y_j) \prod_{i \in S_k \setminus \{j\}} \frac{1}{2}(1 + \alpha_iy_i) \\
    &= \frac{1}{2}(u^{(k)}_\beta + \alpha_j u^{(k)}_\beta y_j) \prod_{i \in S_k \setminus \{j\}} \frac{1}{2}(1 + \alpha_iy_i) \\
    &= \frac{1}{2}(u^{(k)}_\beta + \alpha_j \beta_j u^{(k)}_\beta) \prod_{i \in S_k \setminus \{j\}} \frac{1}{2}(1 + \alpha_iy_i) \\
    &= \frac{1}{2}(u^{(k)}_\beta - u^{(k)}_\beta) \prod_{i \in S_k \setminus \{j\}} \frac{1}{2}(1 + \alpha_iy_i) \\
    &= 0
\end{align*}
where the penultimate equality comes from the fact that $\alpha_j \ne \beta_j$ implies that $\alpha_j\beta_j = -1$. Using the above we have that
\[u^{(k)}_\alpha = u^{(k)}_\alpha \sum_{\beta \in \pm 1^{S_k}_0} \tilde{u}^{(k)}_\beta = u^{(k)}_\alpha\tilde{u}^{(k)}_\alpha = \sum_{\beta \in \pm 1^{S_k}_0} u^{(k)}_\beta \tilde{u}^{(k)}_\alpha = \tilde{u}^{(k)}_\alpha.\]
This completes the proof that $\varphi_1 \circ \varphi_2(u^{(k)}_\alpha) = u^{(k)}_\alpha$ and thus that $\varphi_1 \circ \varphi_2$ is the identity on $C(\qut(G))$. Combining with the above, this proves that $\varphi_1: C^*(\Gamma) \to C(\qut(G))$ and $\varphi_2: C(\qut(G)) \to C^*(\Gamma)$ are isomorphisms of these $C^*$-algebras which are inverse to each other.\\

\noindent\emph{Step 4: Showing that $\Delta_G \circ \varphi_1 = (\varphi_1 \otimes \varphi_1) \circ \Delta_\Gamma$.}\\
Since $\varphi_1$ and the coproducts are $*$-homomorphisms, it suffices to prove the identity on the generators, i.e., that $\Delta_G(\varphi_1(x_i)) = (\varphi_1 \otimes \varphi_1)(\Delta_\Gamma(x_i))$ for all $i \in [n]$. From the above, we have that $\varphi_1(x_i) = y_i = \sum_{\delta \in \pm 1^{S_k}_0} \delta_i u_\delta^{(k)}$ for any $k \in [m]$ such that $i \in S_k$. Thus $(\varphi_1 \otimes \varphi_1)(\Delta_\Gamma(x_i)) = (\varphi_1 \otimes \varphi_1)(x_i \otimes x_i) = y_i \otimes y_i$. We aim to show that $\Delta_G(\varphi_1(x_i)) = \Delta_G(y_i)$ is equal to $y_i \otimes y_i$.

First, let us determine $\Delta_G(u^{(k)}_\delta)$ for $\delta \in \pm 1^{S_k}_0$. Picking $\alpha \in \pm 1^{S_k}_{b_k}$ arbitrarily, we have that $u^{(k)}_\delta = u^{(k)}_{\alpha,\delta \triangle \alpha}$ and thus
\begin{align*}
    \Delta_G\left(u^{(k)}_\delta\right) &= \Delta_G\left(u^{(k)}_{\alpha,\delta \triangle \alpha}\right) \\
    &= \Delta_G\left(u_{(k,\alpha),(k,\delta \triangle \alpha)}\right) \\
    &= \sum_{\substack{l \in [m] \\ \beta \in \pm 1^{S_l}_{b_l}}} u_{(k,\alpha),(l,\beta)} \otimes u_{(l,\beta),(k,\delta \triangle \alpha)} \\
    &= \sum_{\beta \in \pm 1^{S_k}_{b_k}} u_{(k,\alpha),(k,\beta)} \otimes u_{(k,\beta),(k,\delta \triangle \alpha)} \\
    &= \sum_{\delta' \in \pm 1^{S_k}_0} u^{(k)}_{\delta'} \otimes u^{(k)}_{\delta\triangle \delta'}.
\end{align*}
Now we will show that $\Delta_G(y_i) = y_i \otimes y_i$. We have that
\allowdisplaybreaks
\begin{align*}
    \Delta_G(y_i) &= \Delta_G\left(\sum_{\delta \in \pm 1^{S_k}_0} \delta_i u_\delta^{(k)}\right) \\
    &= \sum_{\delta \in \pm 1^{S_k}_0} \delta_i\Delta_G\left(u_\delta^{(k)}\right) \\
    &= \sum_{\delta, \delta' \in \pm 1^{S_k}_0} \delta_i\left( u_{\delta'}^{(k)} \otimes u^{(k)}_{\delta \triangle \delta'}\right) \\
    &= \sum_{\delta, \delta' \in \pm 1^{S_k}_0} \delta_i\delta_i'^2\left( u_{\delta'}^{(k)} \otimes u^{(k)}_{\delta \triangle \delta'}\right) \\
    &= \sum_{\delta, \delta' \in \pm 1^{S_k}_0} \left(\delta'_i u_{\delta'}^{(k)}\right) \otimes \left(\delta_i\delta'_i u^{(k)}_{\delta \triangle \delta'}\right) \\
    &= \sum_{\delta' \in \pm 1^{S_k}_0} \left(\delta'_i u_{\delta'}^{(k)} \otimes \left(\sum_{\delta \in \pm 1^{S_k}_0} (\delta \triangle \delta')_i u_{\delta \triangle \delta'}^{(k)}\right)\right) \\
    &= \sum_{\delta' \in \pm 1^{S_k}_0} \left(\delta'_i u_{\delta'}^{(k)} \otimes \left(\sum_{\delta'' \in \pm 1^{S_k}_0} \delta_i'' u_{\delta''}^{(k)}\right)\right) \\
    &= \sum_{\delta' \in \pm 1^{S_k}_0} \delta'_i u_{\delta'}^{(k)} \otimes y_i \\
    &= y_i \otimes y_i.
\end{align*}
Since we have already shown that $(\varphi_1 \otimes \varphi_1) \circ \Delta_\Gamma(x_i) = y_i \otimes y_i$, this completes the proof.\qeds

\section{Constructing uncolored graphs that have the same quantum automorphism group as a given colored graph}\label{sec:decolor}

In the following, we will describe a procedure that given certain vertex-- and edge-colored graphs produces a decolored graph with isomorphic quantum automorphism group. This procedure is divided into two steps: We first decolor the vertices and then the edges. We start with a lemma known from \cite[Lemma 3.2.3]{Ful}. The degree $\deg(v)$ of a vertex $v$ denotes the number of neighbors of $v$ in the graph $G$.

\begin{lemma}\label{lemdegree}
Let $G$ be a finite graph, $A_G$ be its adjacency matrix and $u_{vw}$, $1\le v,w \le n$ be the generators of $C(\qut(G))$. If $(A_G^l)_{vv}\neq (A_G^l)_{ww}$ for some $l\in \N$, then $u_{vw}=0$. Particularly, if $\deg(v)\neq \deg(w)$, then $u_{vw}=0$.
\end{lemma}

The following lemma can be found in \cite[Lemma 3.2]{disttrans}. The distance $d(v,w)$ between two vertices $v,w \in V(G)$ is the length of a shortest path connecting $v$ and $w$.

\begin{lemma}\label{distance}
Let $G$ be a finite graph and $u_{vw}$, $1\le v,w \le n$ be the generators of $C(\qut(G))$. If we have $d(v,p)\neq d(w,q)$, then $u_{vw}u_{pq}=0$.
\end{lemma}

\begin{lemma}\label{lemdegneighbor}
Let $G$ be a finite graph and $u$ the fundamental representation of $\qut(G)$. Let $v,w \in V(G)$. If there exists $p \in V(G)$ with $d(v,p)=k$ such that $\deg(p) \neq \deg(q)$ for all $q \in V(G)$ with $d(w,q)=k$, then $u_{vw}=0$. 
\end{lemma}

\proof
Let $v,w \in V(G)$ and $p\in V(G)$, $d(v,p)=k$ as above. Using Relation \eqref{rel2} and Lemma \ref{distance}, we get
\begin{align*}
    u_{vw}= u_{vw} \left(\sum_q u_{pq}\right) =u_{vw} \left(\sum_{q;d(w,q)=k} u_{pq}\right).
\end{align*}
We conclude
\begin{align*}
u_{vw}= u_{vw} \left(\sum_{q;d(w,q)=k} u_{pq}\right)=0, 
\end{align*}
since we have $\deg(p) \neq \deg(q)$ and thus $u_{pq}=0$ by Lemma \ref{lemdegree}.
\qed

\begin{remark}
Lemma \ref{lemdegree}, Lemma \ref{distance} and Lemma \ref{lemdegneighbor} also work for vertex-- and edge-colored graphs (take the decolored adjacency matrix in Lemma \ref{lemdegree}), since colors just add more relations on the generators of $C(\qut(G))$.
\end{remark}

We will now define a edge-colored graph $G'$ from the vertex-- and edge-colored graph $G$. 

\begin{definition}\label{defG'}
Let $G$ be a vertex -- and edge-colored graph. Depending on the color $c$, attach a path of length $n_c \in \N_0$ to every vertex colored $c$, where $n_{c_1} \neq n_{c_2}$ for colors $c_1\neq c_2$ and then decolor the vertices of the graph. We choose one of the edge-colors of $G$ and let the edges in the paths all have this edge-color. We denote this new edge-colored (but not vertex-colored) graph by $G'$.
\end{definition}

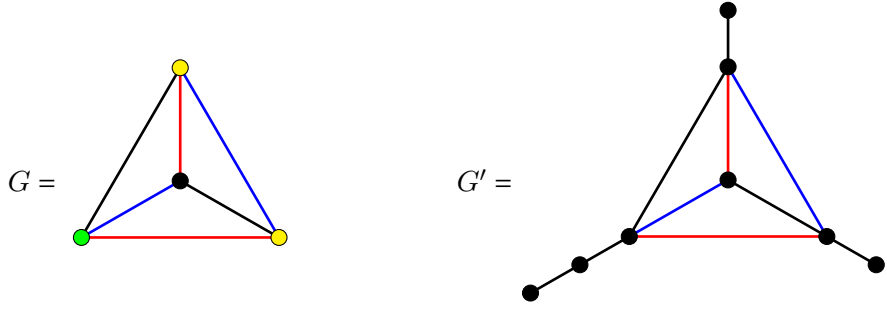
\begin{figure}[h]
\begin{center}\begin{tikzpicture}[scale=1.5, baseline=-1.6cm]
\coordinate[label=left:$G\text{ =}$] (1) at (180:1cm);
\draw[line width=1pt] (330:1cm) -- (0:0cm);
\draw[line width=1pt] (90:1cm) -- (210:1cm);
\draw[blue, line width=1pt] (90:1cm) -- (330:1cm);
\draw[blue, line width=1pt] (210:1cm) -- (0:0cm);
\draw[red, line width=1pt] (90:1cm) -- (0:0cm);
\draw[red, line width=1pt] (210:1cm) -- (330:1cm);
\draw[black,fill=green] (210:1cm) circle (2pt);
\draw[black,fill=black] (0:0cm) circle (2pt);
\foreach \x in {90,330}
{\draw[black,fill=yellow] (\x:1cm) circle (2pt);}
%\draw[black,fill=black] (225:2cm) circle (2pt);
\end{tikzpicture}\hspace{2cm}
\begin{tikzpicture}[scale=1.5]
\coordinate[label=left:$G'\text{ =}$] (1) at (180:1.8cm);
\draw[line width=1pt] (330:1cm) -- (0:0cm);
\draw[line width=1pt] (90:1cm) -- (210:1cm);
\draw[blue, line width=1pt] (90:1cm) -- (330:1cm);
\draw[blue, line width=1pt] (210:1cm) -- (0:0cm);
\draw[red, line width=1pt] (90:1cm) -- (0:0cm);
\draw[red, line width=1pt] (210:1cm) -- (330:1cm);
\draw[black,fill=black] (0:0cm) circle (2pt);
\foreach \x in {90,210,330}
{\draw[black,fill=black] (\x:1cm) circle (2pt);}
\draw[black,fill=black] (210:2cm) circle (2pt);
\draw[black,fill=black] (210:1.5cm) circle (2pt);
\draw[black,fill=black] (90:1.5cm) circle (2pt);
\draw[black,fill=black] (330:1.5cm) circle (2pt);
\draw[line width=1pt] (210:1cm) -- (210:2cm);
\draw[line width=1pt] (90:1cm) -- (90:1.5cm);
\draw[line width=1pt] (330:1cm) -- (330:1.5cm);
\end{tikzpicture}\end{center}\caption{Example for the construction of the graph $G'$ from $G$. We chose the edge-colors between the added paths to be black.}
\end{figure}

The next lemma gives a description of the fundamental representation of $\qut(G')$. We will see in the following theorem that the quantum automorphism groups of $G$ and $G'$ are isomorphic. 

\begin{lemma}\label{lemvertexmagicunitary}
Let $G$ be a vertex -- and edge-colored graph with $\deg(v)\geq 3$ for all $v\in V(G)$ and let $G'$ be the graph as in Definition \ref{defG'}. Denote by $v_i$, $0 \le i \le n_{c(v)}$, the vertices in the added path with $d(v,v_i)=i$ (thus $v=v_0$) and let $V_i=\{x \in V(G'); x=v_i \text{ for some } v\in V(G)\}$. Then the fundamental representation $u$ of $\qut(G')$ is of the form
\begin{align}\label{fundrep1}
\begin{blockarray}{cccccc}
&V_0&V_1&\dots&\dots&V_m\\
\begin{block}{c(ccccc)}
V_0&u_0&0&0&\dots&0\\ V_1&0&u_1&0&\dots&0\\ \vdots&0&0&u_2&\dots&0\\\vdots&\vdots &\vdots&\vdots&\ddots&0\\V_m&0&0&0&0&u_m \\
\end{block}
\end{blockarray}
\end{align}
where additionally $u_{v_iw_i}=u_{v_0w_0}$ for all $i$ and $u_{v_iw_i}=0$ for $c(v)\neq c(w)$, $c(.)$ being the vertex-colors in the original graph $G$.
\end{lemma}

\proof
%Let $u$ be the fundamental representation of $\qut(G')$ and denote by $v_i$, $0 \le i \le n_{c(v)}$, the vertex in the added path with $d(v,v_i)=i$.\dnote{Could this be removed? The notation is already introduced in the Lemma statement.}\\

\noindent\emph{Step 1: It holds $u_{v_iw_0}=u_{v_0w_i}=0$ for $i \neq 0$.}\\
We know that it holds $\deg(v_i) \in \{1,2\}$ for vertices $v_i$ with $i> 0$. Since we have $\deg(w_0)\geq 3$ by assumption, we get $u_{v_iw_0}=0$ by Lemma \ref{lemdegree}. We similarly obtain $u_{v_0w_i}=0$.\\

\noindent\emph{Step 2: We have $u_{v_iw_j}=0$ for $i \neq j$, $i,j >0$.}\\
First assume $i < j$. Then $v_0$ is a vertex with $d(v_i,v_0)=i$ and $\deg(v_0)\geq 3$. Since $i<j$, we know that the vertices $q$ with $d(w_j, q)=i$ are in the path added to $w$ and thus $\deg(q) \in \{1,2\}$. We deduce $u_{v_iw_j}=0$ by Lemma \ref{lemdegneighbor}. The case $i>j$ follows similarly.\\
 
\noindent\emph{Step 3: It holds $u_{v_iw_i}=u_{v_0w_0}$ for all $i$.}\\
We first show $u_{v_0w_0}=u_{v_1w_1}$. It holds
\begin{align*}\allowdisplaybreaks
u_{v_{0}w_{0}}&=u_{v_{0}w_{0}}\left(\sum_{k;(k,w_{0})\in E} u_{v_{1}k}\right)\\
&=u_{v_{0}w_{0}}\left(u_{v_1w_1}+ \sum_{\substack{p_0;(p_0,w_{0})\in E,\\p_0 \in V_0}} u_{v_{1}p_0}\right)\\
&=u_{v_{0}w_{0}}u_{v_1w_1}\\
&=(u_{v_0w_0}+ u_{v_0w_2})u_{v_{1}w_{1}}\\
&=\left(\sum_{k;(k,w_{1})\in E} u_{v_{0}k}\right)u_{v_{1}w_{1}}\\
&=u_{v_1w_1},
\end{align*}
since $u_{v_{1}p_0}=0$ and $u_{v_0w_2}=0$ by \emph{Step 1}. Furthermore, for $i \geq 1$, we have 
\begin{align*}\allowdisplaybreaks
u_{v_{i}w_{i}}&= u_{v_{i}w_{i}}\left(\sum_{k;(k,w_{i})\in E} u_{v_{i+1}k}\right)\\
&=u_{v_{i}w_{i}}(u_{v_{i+1}w_{i+1}}+u_{v_{i+1}w_{i-1}})\\
&=u_{v_{i}w_{i}}u_{v_{i+1}w_{i+1}}\\
&=(u_{v_{i}w_{i}}+ u_{v_{i}w_{i+2}})u_{v_{i+1}w_{i+1}}\\
&= \left(\sum_{k;(k,w_{i+1})\in E} u_{v_{i}k}\right)u_{v_{i+1}w_{i+1}}\\
&=u_{v_{i+1}w_{i+1}}
\end{align*}
since it holds that $u_{v_{i+1}w_{i-1}}=0$ and  $u_{v_{i}w_{i+2}}=0$ by \emph{Step 1} or \emph{Step 2}. Note that if there is no $w_{i+2}$, we still get $u_{v_{i}w_{i}}=u_{v_{i+1}w_{i+1}}$ by a similar calculation, since then $w_i$ is the only neighbor of $w_{i+1}$. 
\\

\noindent\emph{Step 4: It holds $u_{v_iw_i}=0$ for $c(v)\neq c(w)$.}\\
Since $c(v) \neq c(w)$, we know $n_{c(v)}\neq n_{c(w)}$. Assume  $n:=n_{c(v)} < n_{c(w)}$. Then we have $u_{v_{n}w_{n}}=0$ by Lemma \ref{lemdegree}, since $\deg(v_{n})= 1\neq 2=\deg(w_{n})$. We deduce $u_{v_{i}w_{i}}=u_{v_{0}w_{0}}=u_{v_{n}w_{n}}=0$ for all $i$ by \emph{Step 3}. The case $n_{c(v)} > n_{c(w)}$ follows similarly.
\qed

\begin{theorem}\label{thmvertexcolor}
Let $G$ be a vertex -- and edge-colored graph with $\deg(v)\geq 3$ for all $v\in V(G)$ and let $G'$ be as in Definition \ref{defG'}. Then there exists a $*$-isomorphism $\varphi:C(\qut(G))\to C(\qut(G'))$ such that $\Delta_{G'}\circ \varphi=(\varphi\otimes\varphi)\circ \Delta_G$. 
\end{theorem}

\proof
Let $u$ be the fundamental representation of $\qut(G')$. As in Lemma \ref{lemvertexmagicunitary}, we denote by $v_i$, $0 \le i \le n_{c(v)}$ the vertices in the added path with $d(v,v_i)=i$ and let $V_i=\{x \in V(G'); x=v_i \text{ for some } v\in V(G)\}$. Let $c_e$ be the color of the edges in the attached paths in $G'$. For colors $c\neq c_e$, we get that $A_{G'_c}$ is equal to
\begin{align}\label{adjmatrixedge}
\begin{blockarray}{cccccc}
&V_0&V_1&V_2&\dots&V_m\\
\begin{block}{c(ccccc)}
V_0&A_{G_c}&0&0&\dots&0\\ V_1&0&0&0&\dots&0\\V_2&0&0&0&\ddots&0\\\vdots&\vdots &\vdots&\ddots&\ddots&0\\V_m&0&0&0&0& 0\\
\end{block}
\end{blockarray}
\end{align}
By Lemma \ref{lemvertexmagicunitary}, we directly see that $uA_{G_{c}'}=A_{G_{c}'}u$ implies $u_0A_{G_c}=A_{G_c}u_0$ for $c \neq c_e$, $u_0$ being the $(V_0\times V_0)$-block in $u$ (see \eqref{fundrep1}). The adjacency matrix $A_{G_{c_e}'}$ is of the form
\begin{align}\label{adjmatrix1}
\begin{blockarray}{cccccc}
&V_0&V_1&V_2&\dots&V_m\\
\begin{block}{c(ccccc)}
V_0&A_{G_{c_e}}&B_1&0&\dots&0\\ V_1&B_1^t&0&B_2&\dots&0\\V_2&0&B_2^t&0&\ddots&0\\\vdots&\vdots &\vdots&\ddots&\ddots&B_m\\V_m&0&0&0&B_m^t& 0\\
\end{block}
\end{blockarray}
\end{align}
where $B_i$ is the ($V_{i-1} \times V_i$)-matrix with $(B_i)_{a_{i-1}b_i}= \delta_{ab}.$ Also by Lemma \ref{lemvertexmagicunitary}, we see that $uA_{G_{c_e}'}=A_{G_{c_e}'}u$ implies $u_0A_{G_{c_e}}=A_{G_{c_e}}u_0$. Furthermore, since $(u_i)_{v_iw_i}=(u_0)_{v_0w_0}$ for all $i$, we see that $C(\qut(G'))$ is generated by the entries of $u_0$. We conclude that $C(\qut(G'))$ is generated by a magic unitary $u_0$ that fulfills $u_0A_{G_c}=A_{G_c}u_0$ for all edge-colors $c$ and $(u_0)_{v_0w_0}=0$ for $c(v)\neq c(w)$ (see Lemma \ref{lemvertexmagicunitary}). This yields a surjective $*$-homomorphism $\varphi:C(\qut(G))\to C(\qut(G'))$, $w_{ab} \mapsto (u_0)_{ab}$, where $w$ is the fundamental representation of $\qut(G)$.

For the other direction, take the fundamental representation $w$ of $\qut(G)$ and build the matrix $w'$ as follows
\begin{align}\label{constrmagicunitary}
\begin{blockarray}{cccccc}
&V_0&V_1&\dots&\dots&V_m\\
\begin{block}{c(ccccc)}
V_0&w&0&0&\dots&0\\ V_1&0&w_1&0&\dots&0\\ \vdots&0&0&w_2&\dots&0\\\vdots&\vdots &\vdots&\vdots&\ddots&0\\V_m&0&0&0&0&w_m \\
\end{block}
\end{blockarray}
\end{align}
where we put $(w_i)_{a_ib_i}=w_{a_0b_0}$. Note that $w'$ is a magic unitary if and only if the matrices $w_i$ are magic unitaries. We compute \begin{align*}
\sum_{b;n_c(b)\geq i}(w_i)_{a_ib_i}=\sum_{b;n_c(b)\geq i}w_{ab}=\sum_{b}w_{ab}=1,
\end{align*}
where we used $w_{ab}=0$ for $c(a)\neq c(b)$ and the fact that $w$ is a magic unitary. Similarly, we get $\sum_{a;n_c(a)\geq i}(w_i)_{a_ib_i}=1$ and thus $w'$ is a magic unitary. It remains to show that $w'$ commutes with $A_{G'_c}$ for all edge colors $c$. We deduce from \eqref{adjmatrixedge} and \eqref{constrmagicunitary} that $wA_{G_c}=A_{G_c}w$ implies $w'A_{G_{c}'}=A_{G_{c}'}w'$  for $c \neq c_e$. From \eqref{adjmatrix1} and \eqref{constrmagicunitary}, we see that $w'A_{G_{c_e}'}=A_{G_{c_e}'}w'$ is equivalent to $wA_{G_{c_e}}=A_{G_{c_e}}w$, $w_iB_{i+1}=B_{i+1}w_{i+1}$, $w_{i+1}B_{i+1}^t=B_{i+1}^tw_{i}$. Note that we have 
\begin{align*}
    (w_iB_{i+1})_{a_ib_{i+1}}&=\sum_k (w_i)_{a_ik}(B_{i+1})_{kb_{i+1}}\\
    &=(w_i)_{a_ib_i}\\
    &=(w)_{a_0b_0}\\
    &=(w_{i+1})_{a_{i+1}b_{i+1}}\\
    &=\sum_k (B_{i+1})_{a_ik}(w_{i+1})_{kb_{i+1}}\\
    &=(B_{i+1}w_{i+1})_{a_i b_{i+1}}
\end{align*}
by definition of $w_i$ and $B_i$. We similarly get that $w_{i+1}B_{i+1}^t=B_{i+1}^tw_{i}$ is automatically fulfilled by definition of $w_i$ and $B_i$. Since we have $wA_{G_{c_e}}=A_{G_{c_e}}w$ by assumption, we obtain $w'A_{G_{c_e}'}=A_{G_{c_e}'}w'$.
%\dnote{I'm not sure this is easy enough to see to just say it follows from definition. It is straightforward of course but maybe we should write out the expression at least for one of the two cases}. 
Summarizing, we have $w'A_{G_c'}=A_{G_c'}w'$ for all edge-colors $c$. This yields a surjective $*$-homomorphism $\tilde{\varphi}:C(\qut(G'))\to C(\qut(G))$, $u_{xy} \mapsto w'_{xy}$ which is inverse to $\varphi$. Thus, $\varphi$ is a $*$-isomorphism.

It remains to show that $\Delta_{G'}\circ \varphi=(\varphi\otimes\varphi)\circ \Delta_G$. But this is true because
\begin{align*}
    \Delta_{G'}\circ \varphi(w_{ij})&=\Delta_{G'}((u_0)_{ij})\\
    &=\sum_k (u_0)_{ik}\otimes (u_0)_{kj}\\
    &=(\varphi\otimes\varphi)\left(\sum_k w_{ik}\otimes w_{kj}\right)\\
    &=(\varphi\otimes\varphi)\circ \Delta(w_{ij}).
\end{align*}
\qeds
%We show that $u_i$ being magic unitaries is equivalent to $(u_0)_{ab}=0$ for $c(a)\neq c(b)$, $c(.)$ being the colors in the original graph. For one direction, start with the magic unitary $u_m$, which is associated to the vertices $a_m$ with $c(a)=c_{max}$. 
%It holds 
%\begin{align*}
%1=\sum_{k;c(k)=c_{max}}(u_m)_{a_mk_m}=\sum_{k;c(k)=c_{max}}(u_0)_{ak}.
%\end{align*}
% Furthermore, since $u_0$ is a magic unitary, we have $\sum_{k;k=k_0}(u_0)_{ak}=1$. Therefore, we get $\sum_{k;c(k)\neq c_{max}}(u_0)_{ak}=0$ and since all elements in the sum are positive, we deduce $(u_0)_{ab}=0$ for all $a$, $c(a)=c_{max}$ and $b$, $c(b)\neq c_{max}$ (and similarly vice versa). Now, let $u_i$ be the magic unitary with $n_{c_{sec}}=i$ for the second biggest color $c_{sec}$. Take $a$ with $c(a)=c_{sec}$. Then
%%1=\sum_{\substack{k;c(k)=c_{max}\text{ or}\\c(k)=c_{sec}}}(u_i)_{a_ik_i}=\sum_{\substack{k;c(k)=c_{max}\text{ or}\\c(k)=c_{sec}}}(u_0)_{ak}.
%\end{align*}
%Since we have $\sum_{k;k=k_0}(u_0)_{ak}=1$ and $(u_0)_{ab}=0$ for $b$, $c(b)= c_{max}$ by the previous step, we obtain $\sum_{k;c(k)\neq c_{sec}}(u_0)_{ak}=0$. Since all elements are positive, we get $(u_0)_{ab}=0$ for all $a$, $c(a)=c_{sec}$ and $b$, $c(b)\neq c_{sec}$ (and similarly vice versa). Repeating this step for all colors yields $(u_0)_{ab}=0$ for $c(a)\neq c(b)$. 

%\dnote{The above does not handle the colors on the edges of $G$ at all. Since $G$ is colored, it actually has many adjacency matrices (or you can make one matrix with a different value taken in the entries for each color). This should be easy to fix though.}

We will now deal with the edge-colors of $G'$. First, we need the following easy lemma. 

\begin{lemma}\label{lemprojection}
Let $w_{ij}$, $1\le i,j\le n$ be elements in a $C^*$-algebra such that the matrix $w=(w_{ij})_{1\le i,j\le n}$ is a magic unitary. Then $w_{ij}w_{kl}+w_{il}w_{kj}$ is a projection if and only if $w_{ij}w_{kl}=w_{kl}w_{ij}$ and $w_{il}w_{kj}=w_{kj}w_{il}$.
\end{lemma}

\proof
Let $w_{ij}w_{kl}+w_{il}w_{kj}$ be a projection. Since it is self-adjoint, we have 
\begin{align*}
    w_{ij}w_{kl}+w_{il}w_{kj}=w_{kl}w_{ij}+w_{kj}w_{il}.
\end{align*}
Multiplying by $w_{ij}$ and $w_{il}$, respectively, yields $w_{ij}w_{kl}=w_{ij}w_{kl}w_{ij}$ and $w_{il}w_{kj}=w_{il}w_{kj}w_{il}$. But, by taking adjoints, this implies $w_{ij}w_{kl}=w_{kl}w_{ij}$ and $w_{il}w_{kj}=w_{kj}w_{il}$. The other direction is clear, since $w_{ij}w_{kl}+w_{il}w_{kj}$ is the sum of two orthogonal projections if $w_{ij}w_{kl}=w_{kl}w_{ij}$ and $w_{il}w_{kj}=w_{kj}w_{il}$. 
\qeds

%\dnote{Perhaps there should be some remark about how $u_{e_i f_i}$ in the lemma below does not depend on the order of $(v,w)$ or $(x,y)$ where $e = (v,w)$ and $f = (x,y)$ since swapping the order does not change the resulting operator.}

We will define a (uncolored) graph $G''$ from the edge-colored graph $G'$. In the next theorem, we will then see that for certain graphs, the quantum automorphism groups of $G'$ and $G''$ are isomorphic. Let $G$ be a graph and $e=(u,v)\in E(G)$. We say that we subdivide $e$ if we delete the edge $e=(u,v)$ from $G$ and add a vertex $w$ as well as edges $(u,w)$ and $(w,v)$ to the graph.

%\dnote{Is $c_0$ in the following the color used for the edges added to $G$ to obtain $G'$, it is a little unclear.}

\begin{definition}\label{defG''}
Let $G$ be a vertex -- and edge-colored graph. First decolor the vertices by applying the construction as in Definition \ref{defG'} to obtain the edge-colored graph $G'$. We denote the color of the newly added edges in the graph $G'$ by $c_0$. We subdivide each colored edge with $c(e)\neq c_0$ and add a path of length $m_c$ to the subdivision, where $m_{c_1} \neq m_{c_2}$ for colors $c_1\neq c_2$  and then decolor the edges in the graph $G'$. We call this graph $G''$.
\end{definition}

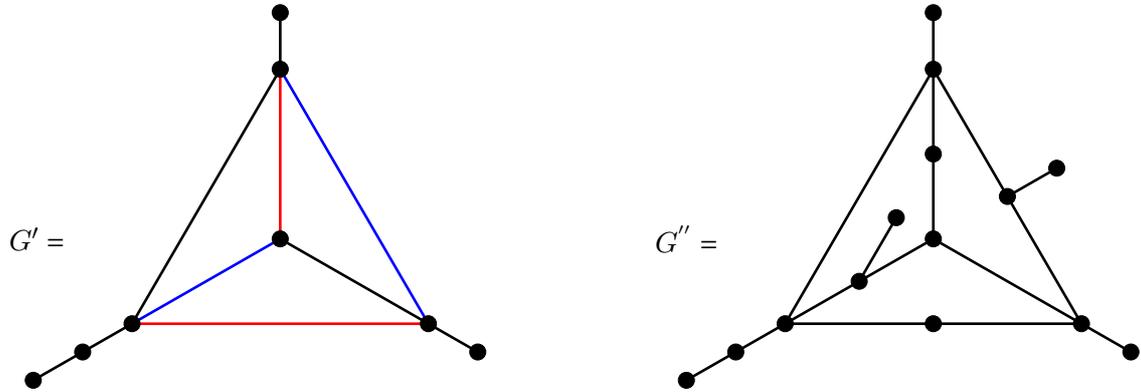
\begin{figure}[h]
\begin{center}
\begin{tikzpicture}[scale=1.5]
\coordinate[label=left:$G'\text{ =}$] (1) at (180:1.8cm);
\draw[line width=1pt] (330:1.5cm) -- (0:0cm);
\draw[line width=1pt] (90:1.5cm) -- (210:1.5cm);
\draw[blue, line width=1pt] (90:1.5cm) -- (330:1.5cm);
\draw[blue, line width=1pt] (210:1.5cm) -- (0:0cm);
\draw[red, line width=1pt] (90:1.5cm) -- (0:0cm);
\draw[red, line width=1pt] (210:1.5cm) -- (330:1.5cm);
\draw[black,fill=black] (0:0cm) circle (2pt);
\foreach \x in {90,210,330}
{\draw[black,fill=black] (\x:1.5cm) circle (2pt);}
\draw[black,fill=black] (210:2.5cm) circle (2pt);
\draw[black,fill=black] (210:2cm) circle (2pt);
\draw[black,fill=black] (90:2cm) circle (2pt);
\draw[black,fill=black] (330:2cm) circle (2pt);
\draw[line width=1pt] (210:1.5cm) -- (210:2.5cm);
\draw[line width=1pt] (90:2cm) -- (90:1.5cm);
\draw[line width=1pt] (330:2cm) -- (330:1.5cm);
\end{tikzpicture}\hspace{2cm}
\begin{tikzpicture}[scale=1.5] %baseline=-1.1cm]
\coordinate[label=left:$G^{''}\text{ =}$] (1) at (180:1.8cm);
\draw[line width=1pt] (330:1.5cm) -- (0:0cm);
\draw[line width=1pt] (90:1.5cm) -- (210:1.5cm);
\draw[line width=1pt] (90:1.5cm) -- (330:1.5cm);
\draw[line width=1pt] (210:1.5cm) -- (0:0cm);
\draw[line width=1pt] (90:1.5cm) -- (0:0cm);
\draw[line width=1pt] (210:1.5cm) -- (330:1.5cm);
\draw[black,fill=black] (0:0cm) circle (2pt);
\foreach \x in {90,210,330}
{\draw[black,fill=black] (\x:1.5cm) circle (2pt);}
\draw[black,fill=black] (210:2.5cm) circle (2pt);
\draw[black,fill=black] (210:2cm) circle (2pt);
\draw[black,fill=black] (90:2cm) circle (2pt);
\draw[black,fill=black] (330:2cm) circle (2pt);
\draw[line width=1pt] (210:1.5cm) -- (210:2.5cm);
\draw[line width=1pt] (90:2cm) -- (90:1.5cm);
\draw[line width=1pt] (330:2cm) -- (330:1.5cm);
%\draw[black,fill=black] (225:2cm) circle (2pt);
\draw[black,fill=black] (90:0.75cm) circle (2pt);
\draw[black,fill=black] (270:0.75cm) circle (2pt);
\draw[black,fill=black] (210:0.75cm) circle (2pt);
\draw[black,fill=black] (30:0.75cm) circle (2pt);
\draw[black,fill=black] (30:1.25cm) circle (2pt);
\draw[line width=1pt] (30:0.75cm) -- (30:1.25cm);
\draw[line width=1pt] (210:0.75cm) -- (150:0.375cm);
\draw[black,fill=black] (150:0.375cm) circle (2pt);
\end{tikzpicture}
\end{center}\caption{Example for the construction of the graph $G^{''}$ from $G'$. We chose black to be the edge-color $c_0$.}
\end{figure}

\begin{lemma}\label{lemedgemagicunitary}
Let $G$ be a vertex -- and edge-colored graph such that $\deg(v)\geq 3$ for all $v\in V(G)$ and let $G''$ be the graph as in Definition \ref{defG''}. We denote the vertex that subdivided the edge $e$ in $G$ by $e_0$ and the vertices in the added path with $d(e_0,e_i)=i$ by $e_i$, $1 \le i \le m_{c(e)}$. Furthermore $E_i'=\{x \in V(G''); x=e_i \text{ for some } e\in E(G)\}$. Then, the fundamental representation $u$ of $\qut(G'')$ is of the form
\begin{align}\label{magicunitary2}
\begin{blockarray}{cccccc}
&V(G')&E_0'&\dots&\dots&E_m'\\
\begin{block}{c(ccccc)}
V(G')&w&0&0&\dots&0\\ E_0'&0&u_0&0&\dots&0\\ \vdots&0&0&u_1&\dots&0\\\vdots&\vdots &\vdots&\vdots&\ddots&0\\E_m'&0&0&0&0&u_m \\
\end{block}
\end{blockarray}
\end{align}
where $u_{e_if_i}=u_{vx}u_{wy}+u_{vy}u_{wx}$ for $c(e)=c(f)$, $e=(v,w), f=(x,y)$ and $u_{e_if_i}=0$ for $c(e)\neq c(f)$.
\end{lemma}

\proof
%Let $u$ be the fundamental representation of $\qut(G'')$ and denote the vertex that subdivided the edge $e$ in $G$ by $e_0$ and the vertices in the added path with $d(e_0,e_i)=i$ by $e_i$, $1 \le i \le n_{c(e)}$\dnote{This is already stated in the Lemma statement.}. Moreover, we denote by $v_i$, $0 \le i \le n_{c(v)}$ the vertices in the added path with $d(v,v_i)=i$.\\%note that subdividing does not change the degree of vertices

\noindent\emph{Step 1: It holds $u_{e_i v_j}=u_{v_j e_i}=0$ for $v_j\in V(G')$, $e_i \in E_i'$.}\\
First assume $i>0, j=0$. Then $\deg(e_i)\in \{1,2\}$ and $\deg(v_0)\geq 3$, therefore $u_{e_iv_0}=0$ by Lemma \ref{lemdegree}.

Now, let $i=0, j>0$. If $m_{c(e_0)}>0$, then $\deg(e_0)=3\neq \deg(v_j)$, since $\deg(v_j)\in \{1,2\}$. Therefore $u_{e_0v_j}=0$ by Lemma \ref{lemdegree}. Let $m_{c(e_0)}=0$. If $n_{c(v_0)}=1$, then $\deg(v_1)=1 \neq 2=\deg(e_0)$ and thus $u_{e_0v_1}=0$ by Lemma \ref{lemdegree}. If $n_{c(v_0)}\geq 2$, then $v_j$ has at least one neighbor of degree $1$ or $2$. Since $e_0$ only has neighbors of degree $\geq 3$, we get $u_{e_0v_j}=0$ by Lemma \ref{lemdegneighbor}.

Assume $i\neq j$ and $i,j \neq 0$. Let furthermore $i<j$. Then $e_0$ is a vertex with $d(e_i,e_0)=i$ and $\deg(e_0)= 3$. Since $i<j$, we know that the vertices $q$ with $d(v_j, q)=i$ are in the path and thus $\deg(q) \in \{1,2\}$. We deduce $u_{e_iv_j}=0$ by Lemma \ref{lemdegneighbor}. The case $i>j$ follows similarly.

It remains to show $u_{e_iv_i}=0$. If $n_{c(v_0)}>0$, then $\deg(v_0)\geq 4$ and $d(v_0,v_i)=i$. We know $\deg(q)\le 3$ for vertices $q$ with $d(q,e_i)=i$, since either $q=e_0$ or $q$ is a vertex in the added path. We deduce $u_{e_iv_i}=0$ by Lemma \ref{lemdegneighbor}. Now assume $n_{c(v_0)}=0$. If $m_{c(e_0)}=0$, then we know $\deg(e_0)=2$ and thus $u_{e_0v_0}=0$ by Lemma \ref{lemdegree}. If $m_{c(e_0)}>0$, then $e_0$ has a neighbor of degree $1$ or $2$. If $v_0$ has no neighbor of this degree, then $u_{e_0v_0}=0$ by Lemma \ref{lemdegneighbor}. The vertex $v_0$ only has a neighbor of degree $2$ if there exists a subdivision $f_0$ of some edge $f = (v_0,w_0)$ with $m_{c(f)}=0$. Then, it holds
\begin{align*}
u_{e_0v_0}= u_{e_0v_0}\left(\sum_{k;(k,e_0)\in E(G'')}u_{kf_0}\right)=u_{e_0v_0}u_{e_1f_0}.
\end{align*}
If $m_{c(e_0)}=1$, then $\deg(e_1)=1 \neq 2 = \deg(f_0)$ which yields $u_{e_1f_0}=0$ by Lemma \ref{lemdegree}. If $m_{c(e_0)}\geq 2$, then $e_1$ has a neighbor $e_2$ with $\deg(e_2) \in \{1,2\}$. Since the neighbors $v_0$, $w_0$ of $f_0$ have degree $\geq 3$, we get $u_{e_1f_0}=0$ by Lemma \ref{lemdegneighbor}. We deduce $u_{e_0v_0}=0$ in all those cases. 
\\

\noindent\emph{Step 2: It holds $u_{e_i f_0}=u_{e_0f_i}=0$ for $i \neq 0$.}\\
We know that it holds $\deg(e_i) \in \{1,2\}$ for vertices $e_i$ with $i> 0$. If $m_{c(f_0)}>0$, then we have $\deg(f_0)= 3$ and get $u_{e_if_0}=0$ by Lemma \ref{lemdegree}. Let $m_{c(f_0)}=0$. If $m_{c(e_0)}=1$, then $\deg(e_1)=1$ and thus $u_{e_1f_0}=0$ by Lemma \ref{lemdegree}. If $m_{c(e_0)}\geq 2$, then $e_i$ has a neighbor with degree $1$ or $2$. Since the neighbors of $f_0$ have degree $\geq 3$ by assumption, we get $u_{e_if_0}=0$ by Lemma \ref{lemdegneighbor}. We similarly obtain $u_{e_0f_i}=0$.\\

\noindent\emph{Step 3: It holds $u_{e_if_j}=0$ for $i \neq j$, $i,j>0$.}\\
 First assume $i < j$. Then $e_0$ is a vertex with $d(e_i,e_0)=i$ and $\deg(e_0)= 3$. Since $i<j$, we know that the vertices $q$ with $d(f_j, q)=i$ are in the path and thus $\deg(q) \in \{1,2\}$. We deduce $u_{e_if_j}=0$ by Lemma \ref{lemdegneighbor}. The case $i>j$ follows similarly.\\
 
 \noindent\emph{Step 4: It holds  $u_{e_if_i}=u_{e_0f_0}$ for all $i$.}\\
 We first show $u_{e_0f_0}=u_{e_1f_1}$. It holds 
\begin{align*}
u_{e_{0}f_{0}}&=u_{e_{0}f_{0}}\left(\sum_{k;(k,f_{0})\in E(G'')} u_{e_{1}k}\right)\\
&=u_{e_{0}f_{0}}\left(u_{e_1f_1}+ \sum_{p_0;(p_0,f_{0})\in E(G''), p_0 \in V(G')} u_{e_{1}p_0}\right)\\
&=u_{e_{0}f_{0}}u_{e_1f_1}\\
&=(u_{e_0f_0}+ u_{e_0f_2})u_{e_{1}f_{1}}\\
&=\left(\sum_{k;(k,f_{1})\in E(G'')} u_{e_{0}k}\right)u_{e_{1}f_{1}}\\
&=u_{e_1f_1},
\end{align*}
since $u_{e_{1}p_0}=0$ by \emph{Step 1} and $u_{e_0f_2}=0$ by \emph{Step 2}. Furthermore, for $i \geq 1$, we have 
\begin{align*}
u_{e_{i}f_{i}}&= u_{e_{i}f_{i}}\left(\sum_{k;(k,f_{i})\in E(G'')} u_{e_{i+1}k}\right)\\
&=u_{e_{i}f_{i}}(u_{e_{i+1}f_{i+1}}+u_{e_{i+1}f_{i-1}})\\
&=u_{e_{i}f_{i}}u_{e_{i+1}f_{i+1}}\\
&=(u_{e_{i}f_{i}}+ u_{e_{i}f_{i+2}})u_{e_{i+1}f_{i+1}}\\
&= \left(\sum_{k;(k,f_{i+1})\in E(G'')} u_{e_{i}k}\right)u_{e_{i+1}f_{i+1}}\\
&=u_{e_{i+1}f_{i+1}}
\end{align*}
since we have $u_{e_{i+1}f_{i-1}}=0$ and  $u_{e_{i}f_{i+2}}=0$ by \emph{Step 2} or \emph{Step 3}. Note that if there is no $f_{i+2}$, we still get $u_{e_{i}f_{i}}=u_{e_{i+1}f_{i+1}}$ by a similar calculation, since then $f_i$ is the only neighbor of $f_{i+1}$.\\

\noindent\emph{Step 5: It holds $u_{e_if_i}=0$ for $c(e)\neq c(f)$.}\\
Since $c(e) \neq c(f)$, we know $m_{c(e)}\neq m_{c(f)}$. First assume  $m:=m_{c(e)} < m_{c(f)}$. Then we have $u_{e_{m}f_{m}}=0$ by Lemma \ref{lemdegree}, since $\deg(e_{m})= 1\neq 2=\deg(f_{m})$. We deduce $u_{e_{i}f_{i}}=u_{e_{0}f_{0}}=u_{e_{m}f_{m}}=0$ for all $i$ by \emph{Step 3}. The case $m_{c(e)} > m_{c(f)}$ is similar.
\\

\noindent\emph{Step 6: It holds $u_{e_if_i}=u_{vx}u_{wy}+u_{vy}u_{wx}$ for $e=(v,w)$, $f=(x,y)$.}\\
We have 
\begin{align*}
u_{e_0f_0}&= \left(\sum_{k;(k,f_0)\in E(G'')}u_{vk}\right)u_{e_0f_0}\left(\sum_{l;(l,f_0)\in E(G'')}u_{wl}\right)\\
&=(u_{vx}+u_{vy}+u_{vf_1})u_{e_0f_0}(u_{wx}+u_{wy}+u_{wf_1})\\
&=(u_{vx}+u_{vy})u_{e_0f_0}(u_{wx}+u_{wy}),
\end{align*}
where we used $u_{vf_1}=0$, $u_{wf_1}=0$ by \emph{Step 1}. Note that $e_0$ is the only vertex in $E_0'$ that is a common neighbor of $v$ and $w$. Thus, we have $(u_{vx}+u_{vy})u_{kf_0}(u_{wx}+u_{wy})=0$ for all $k \neq e_0$ by Lemma \ref{distance} and \emph{Step 1}. Therefore, we get 
\begin{align*}
u_{e_0f_0}&=(u_{vx}+u_{vy})u_{e_0f_0}(u_{wx}+u_{wy})\\
&=(u_{vx}+u_{vy})\left(\sum_{k;(k,v)\in E(G''), (k,w)\in E(G'')}u_{kf_0}\right)(u_{wx}+u_{wy})\\
&=(u_{vx}+u_{vy})(u_{wx}+u_{wy})\\
&=u_{vx}u_{wy}+u_{vy}u_{wx}.
\end{align*}
By \emph{Step 4}, we obtain $u_{e_if_i}=u_{vx}u_{wy}+u_{vy}u_{wx}$.
\qeds

\begin{remark}
Note that the operator $u_{(v,w)(x,y)}=u_{vx}u_{wy}+u_{vy}u_{wx}$ does not depend on the order of $(v,w)$ or $(x,y)$, since
\begin{align*}
    u_{(v,w)(x,y)}&=u_{vx}u_{wy}+u_{vy}u_{wx}=u_{vy}u_{wx}+u_{vx}u_{wy}=u_{(v,w),(y,x)},\\
    u_{(v,w)(y,x)}&=u_{vy}u_{wx}+u_{vx}u_{wy}=u_{wx}u_{vy}+u_{wy}u_{vx}=u_{(w,v),(x,y)}.\\
\end{align*}
\end{remark}

Before stating the theorem, we first need to define a quantum subgroup of the quantum automorphism group of the graph $G'$. It is straightforward to check that the comultiplication is a *-homomorphism.  

\begin{definition}\label{defQut*}
Let $G$ be a vertex -- and edge-colored graph and $G'$ as in Definiton \ref{defG'}. We define $\qut^*_{c_0}(G')$ to be the compact matrix quantum group whose corresponding $C^*$-algebra is generated by a magic unitary $x$ with $xA_{G_c}=A_{G_c}x$ for every edge color $c$ and $x_{ik}x_{jl}=x_{jl}x_{ik}$ for $c((i,j))=c((k,l))$, $c\neq c_0$, where $c_0$ is the edge-color we choose for the newly added edges in $G'$.
\end{definition}

\begin{theorem}\label{theoremedgecolor}
Let $G$ be a vertex -- and edge-colored graph such that $\deg(v)\geq 3$ for all $v\in V(G)$. Let $G''$ be the graph as in Definition \ref{defG''} and let $\qut^*_{c_0}(G')$ be the compact matrix quantum group as in Definition \ref{defQut*}. Then there exists a $*$-isomorphism $\varphi:C(\qut^*_{c_0}(G')) \to C(\qut(G''))$ such that $\Delta_{G''}\circ \varphi=(\varphi\otimes\varphi)\circ \Delta_{G'}$.
\end{theorem}
%\dnote{In the statement you use $w$ as the magic unitary whose entries generate $C(\qut^*(G'))$, but you also use $w$ as a block of $u$. So probably should change one of them. Are we taking $c_0$ to be the color of the new edges added to $G$ to obtain $G'$?}
\proof
Let $u$ be the fundamental representation of $\qut(G'')$. As in Lemma \ref{lemedgemagicunitary}, we denote the vertex that subdivided the edge $e$ in $G$ by $e_0$ and the vertices in the added path with $d(e_0,e_i)=i$ by $e_i$, $1 \le i \le m_{c(e)}$. Moreover, we denote by $v_i$, $0 \le i \le n_{c(v)}$ the vertices in the added path with $d(v,v_i)=i$. The adjacency matrix of $G''$ is of the form
\begin{align}\label{adjmatrix2}
\begin{blockarray}{cccccc}
&V(G')&E_0'&E_1'&\dots&E_m'\\
\begin{block}{c(ccccc)}
V(G')&A_{G'_0}&B_{G''}&0&\dots&0\\ E_0'&B_{G''}^t&0&C_1&\dots&0\\E_1'&0&C_1^t&0&\ddots&0\\\vdots&\vdots &\vdots&\ddots&\ddots&C_m\\E_m'&0&0&0&C_m^t& 0\\
\end{block}
\end{blockarray}
\end{align}
where $A_{G_0}$ is the adjacency matrix of the edge-color that is not subdivided together with the paths from the construction of $G'$, $B_{G''}$ is the matrix with 
\begin{align*}
(B_{G''})_{ve}=\begin{cases} 1 \quad \text{  for $v$ incident to $e$,} \\ 0 \quad \text{  otherwise,}\end{cases}
\end{align*}
 and $C_i$ is the ($E_{i-1}' \times E_i'$)-matrix with $(C_i)_{a_{i-1}b_i}= \delta_{ab}$. %\dnote{Need $c(e) \ne c_0$ in the above I think.} There is no vertex e with color c_0 in E_0
 Using Lemma \ref{lemedgemagicunitary}, we see that  $uA_{G''}=A_{G''}u$ implies $wA_{G'_0}=A_{G'_0}w$, $wB_{G''}=B_{G''}u_0$, $u_0B_{G''}^t=B_{G''}^tw$, where $w$ and $u_0$ are blocks in $u$ (see \eqref{magicunitary2}). Furthermore, since $(u_i)_{e_if_i}=w_{ak}w_{bl}+w_{al}w_{bk}$ for $e=(a,b), f=(k,l)$, we also have that $C(\qut(G''))$ is generated by the entries of the matrix $w$. We will now show that $w$ fulfills the relations of the generators of $C(\qut^*_{c_0}(G'))$.
 We already know $wA_{G'_0}=A_{G'_0}w$. We will now show $wA_{G_c}=A_{G_c}w$ for every edge color $c\neq c_0$. By Lemma \ref{lem:productrelation}, this is equivalent to $w_{ik}w_{jl}=0$ for $c((i,j))\neq c((k,l))$ if $(i,j)\in E_0'$ or $(k,l)\in E_0'$. We will show that the elements of $w$ fulfill those relations. Let $v \in V(G')$, $(a,b)=e \in E_0'$. 
 It holds
 \begin{align}
     w_{va}+w_{vb}&=\sum_{k;k \text{ inc. to }e} w_{vk}\nonumber\\
     &=\sum_k w_{vk}(B_{G''})_{ke}\nonumber\\
     &=(wB_{G''})_{ve}\nonumber\\
     &=(B_{G''}u_0)_{ve}\nonumber\\
     &=\sum_{k}(B_{G''})_{vk}(u_0)_{ke}\nonumber\\
     &=\sum_{k;k \text{ inc. to }v}(u_0)_{ke}\nonumber\\
     &=\sum_{s;(v,s)\in E_0}(w_{va}w_{sb}+w_{vb}w_{sa}).\label{eqst8}
 \end{align}
Multiplying $w_{va}$ from left and right yields
\begin{align*}
    w_{va}=\sum_{s;(v,s)\in E_0} w_{va}w_{sb}w_{va}
\end{align*}
We deduce $\sum_{s;(v,s)\notin E_0'} w_{va}w_{sb}w_{va}=0$ and since all elements in the sum are positive, we get $w_{va}w_{sb}w_{va}=0$ for all $v,s \in V(G')$ with $(v,s) \notin E_0'$. We deduce $w_{va}w_{sb}=0$ for all $v,s \in V$ with $(v,s) \notin E_0'$. Similarly, by using $u_0B_{G''}^t=B_{G''}^tw$, we get $w_{av}w_{bs}=0$ for all $v,s \in V(G')$ with $(v,s) \notin E_0'$. %\dnote{Do I understand correctly that the edges not in $E_0$ are just the edges that were colored $c_0$, and thus were not subdivided? If this is the case then isn't it immediate that $w_{ia}w_{xb} =0$ if $(i,x) \notin E_0$ and $(a,b) \in E_0$ since the distances in $G''$ are different?} Also non-edges are not in E_0, therefore the distance argument does not work for all (i,x) \notin E_0
It remains to show $w_{ik}w_{jl}=0$ for $c(e)\neq c(f)$, $e=(i,j)\in E_0'$, $f=(k,l)\in E_0'$.
We know $(u_0)_{ef}=0$ for $c(e)\neq c(f)$ from Lemma \ref{lemedgemagicunitary}. Then 
\begin{align*}
    w_{ik}w_{jl}+w_{il}w_{jk}=(u_0)_{ef}=0
\end{align*}
and by multiplying $w_{ik}$ and $w_{il}$ from the left, respectively, we get $ w_{ik}w_{jl}=0$ and $w_{il}w_{jk}=0$. Thus, we have shown $w_{ik}w_{jl}=0$ for $c((i,j))\neq c((k,l))$ if $(i,j)\in E_0'$ or $(k,l)\in E_0'$ which is equivalent to $wA_{G_c}=A_{G_c}w$ for every edge color $c \neq c_0$. 

Summarizing, we get $wA_{G_c}=A_{G_c}w$ for every edge color $c$.
By Lemma \ref{lemedgemagicunitary}, we furthermore know that $(u_0)_{ab}=w_{ik}w_{jl}+w_{il}w_{jk}$ for $a = (i,j)$, $b = (k,l)$ and $c(a)=c(b)$, $c(a)\neq c_0$ and thus $w_{ik}w_{jl}+w_{il}w_{jk}$ is a projection. Then Lemma \ref{lemprojection} yields $w_{ij}w_{kl}=w_{kl}w_{ij}$ and $w_{il}w_{kj}=w_{kj}w_{il}$.
Therefore, we get a surjective $*$-homomorphism $\varphi:C(\qut^*_{c_0}(G')) \to C(\qut(G'')$, $w'_{ij}\mapsto w_{ij}$, where $w'$ is the fundamental representation of $\qut^*_{c_0}(G')$.

For the other direction, take the fundamental representation $w'$ of $\qut^*_{c_0}(G')$ and build the matrix $w''$ as follows
\begin{align*}
\begin{blockarray}{cccccc}
&V(G')&E_0'&\dots&\dots&E_m'\\
\begin{block}{c(ccccc)}
V(G')&w'&0&0&\dots&0\\ E_0'&0&u_0'&0&\dots&0\\ \vdots&0&0&u'_1&\dots&0\\\vdots&\vdots &\vdots&\vdots&\ddots&0\\E_m'&0&0&0&0&u'_m \\
\end{block}
\end{blockarray}
\end{align*}
where we put $(u'_i)_{e_if_i}=w'_{ak}w'_{bl}+w'_{al}w'_{bk}$ for $e=(a,b), f=(k,l)$. Those elements are projections by Lemma \ref{lemprojection}, because we have $w'_{ak}w'_{bl}=w'_{bl}w'_{ak}$ for $(a,b),(k,l)\in E_0'$ (for the same color by assumption, for different colors because the product is $0$). Note that $w''$ is a magic unitary if and only if the matrices $u'_i$ are magic unitaries. Let $e=(a,b)\in E_0'$. We compute \begin{align*}
\sum_{f;n_c(f)\geq i}(u'_i)_{e_if_i}=\sum_{k,l;n_c((k,l))\geq i}(w'_{ak}w'_{bl}+w'_{al}w'_{bk})=\sum_{k,l; (k,l) \in E_0'}(w'_{ak}w'_{bl}+w'_{al}w'_{bk})=1,
\end{align*}
where we used $w'_{ak}w'_{bl}=0$ for $c((a,b))\neq c((k,l))$ and the fact that $w'$ is a magic unitary. %\dnote{I think the sums should be over $\{l,k\} \in E_0$, otherwise you get 2 instead of 1}. 
We get $\sum_{f;n_c(f)\geq i}(u'_i)_{f_ie_i}=1$ similarly and therefore $w''$ is a magic unitary. It remains to show that $w''$ commutes with $A_{G''}$.
Similar to equation \eqref{eqst8}, we see that%\dnote{You use $j$ to denote the edge $(a,b)$ in the following I think, but you never specify that. Also we are usually using $j$ as a vertex so it might be better to use something else.}
\begin{align*}
    (w'B_{G''})_{i(a,b)}=w'_{ia}+w'_{ib} \quad \text{ and }\quad  (B_{G''}u'_0)_{i(a,b)}=\sum_{x;(i,x)\in E_0'}(w'_{ia}w'_{xb}+w'_{ib}w'_{xa}).
\end{align*}
Using $w'_{ia}w'_{xb}=0$ for $(i,x)\notin E_0', (a,b)\in E_0'$, we get
\begin{align*}
    \sum_{x;(i,x)\in E_0'}(w'_{ia}w'_{xb}+w'_{ib}w'_{xa})=\sum_{x}(w'_{ia}w'_{xb}+w'_{ib}w'_{xa})=w'_{ia}+w'_{ib},
\end{align*}
and therefore $w'B_{G''}= B_{G''}u'_0$. We obtain $u'_0B_{G''}^t=B_{G''}^tw'$ similarly.
From \eqref{adjmatrix2} and the form of $w''$, we see that $w''A_{G''}=A_{G''}w''$ is equivalent to $w'A_{G'_0}=A_{G'_0}w'$, $w'B_{G''}=B_{G''}u'_0$, $u'_0B_{G''}^t=B_{G''}^tw'$, $u'_iC_{i+1}=C_{i+1}u'_{i+1}$, $u'_{i+1}C_{i+1}^t=C_{i+1}^tu'_{i}$. Note that $u'_i C_{i+1}=C_{i+1}u'_{i+1}$, $u'_{i+1}C_{i+1}^t=C_{i+1}^t u'_{i}$ are automatically fulfilled by definition of $u'_i$ and $C_i$. Therefore, we get $w''A_{G''}=A_{G''}w''$. This yields a surjective $*$-homomorphism $\tilde{\varphi}:C(\qut(G''))\to C(\qut^*(G')_{c_0})$, $u_{ij} \mapsto w''_{ij}$ which is inverse to $\varphi$. Thus, $\varphi$ is a $*$-isomorphism.

It remains to show that $\Delta_{G''}\circ \varphi=(\varphi\otimes\varphi)\circ \Delta_{G'}$, which is true because of the following:
\begin{align*}
    \Delta_{G''}\circ \varphi(w'_{ij})&=\Delta_{G'}(w_{ij})\\
    &=\sum_k w_{ik}\otimes w_{kj}\\
    &=(\varphi\otimes\varphi)\left(\sum_k w'_{ik}\otimes w'_{kj}\right)\\
    &=(\varphi\otimes\varphi)\circ \Delta_{G'}(w'_{ij}).
\end{align*}
\qed\vspace{0.5cm}

We get the following corollary, which we will use in the next section to construct a graph with quantum symmetry and finite quantum automorphism group. 

\begin{cor}\label{coreq}
Let $G$ be a vertex -- and edge-colored graph such that $\deg(v)\geq 3$ for all $v\in V(G)$. Denote by $G'$ and $G''$ the graphs constructed in Definitions \ref{defG'} and \ref{defG''}. If $\qut^*_{c_0}(G')\cong \qut(G')$, then $\qut(G'')\cong \qut(G)$. 
\end{cor}

\section{(Uncolored) graphs whose quantum automorphism group is the dual of a solution group}\label{sec:qsymfinaut}

In this section, we will look at certain graphs from Definition \ref{defgraph}, where we replace one of the edge colors by non-edges. Furthermore, we restrict to graphs coming from linear constraint systems as in the following definition. By using Theorem \ref{thm:qaut2cgamma} and the decoloring procedure in Section \ref{sec:decolor}, we will obtain (decolored) graphs whose quantum automorphism group is the dual of the corresponding solution group. 

%Cite arkhipov?
\begin{definition}\label{ArkhipovLBCS}
Let $H$ be a connected graph with vertex set $[m]$ and label the edges $1,\dots, n:=|E(H)|$. Let $M_H \in \mathbb{F}_2^{m \times n}$ be the matrix, where
\begin{align*}
    (M_H)_{ki}= \begin{cases}1 & \text{if } k\in V(H) \text{ is incident to } i \in E(H);\\ 0 & \text{o.w.}\end{cases}
\end{align*}
%Consider the linear system $M_Hx=0$. We denote the corresponding solution group by $\Gamma_H:=\Gamma(M_H,0).$
\end{definition}

\begin{definition}\label{defnewgraph}
Let $H$ be a connected graph, let $M_H$ be as in Definition \ref{ArkhipovLBCS} and let $G(M_H, b)$ be the colored graph as in Definition \ref{defgraph}. Recall $S_k(M_H)=\{i \in [n] : (M_H)_{ki} = 1\}$. Note that it holds $|S_l(M_H)\cap S_k(M_H)|\le 1$ for all $l\neq k$, since otherwise there would be multiple edges between $l$ and $k$ in the graph $H$ by definition of $M_H$. By Remarks \ref{rem1} and \ref{rem2}, we know that we can color the edges between any $V_l$ and $V_k$, $k\neq l$, by the same $(2^{|S_k\cap S_l|}-1)$ numbers. In our case, this means that for every $|S_l\cap S_k|\neq 0$, we can use the same two numbers for any $V_l$, $V_k$. We color the edge between $(l,\alpha)$, $(k,\beta)$ by $\alpha_i\beta_i \in \{-1,1\}$, where $i\in [n]$ is the number for which $(M_H)_{ki}=(M_H)_{li}=1$. We then define the graph $G_*(M_H,b)$ to be the graph obtained from the previous graph by replacing the edges of color $1$ by non-edges. %\dnote{Strictly speaking $-1$ is not an element of $\pm 1^{S_l \cap S_k}$, rather the function $f:\{a\} \to \{\pm 1\}$, $f(a) = -1$ is, where $a$ is the single element in $S_l \cap S_k$. This means that this color is actually different for different $l$ and $k$. I think we may need to be a bit more clear about this. In particular, we need to make clear that in the following lemma we are taking all edges between any two $V_l$ and $V_k$ with $S_l \cap S_k \ne \varnothing$ to be the same color and also using this color on the added paths.}
\end{definition} 

\begin{remark}
As already mentioned in Remark \ref{rem2}, Theorem \ref{thm:qaut2cgamma} is true if we replace $G(M_H,0)$ by $G_*(M_H,0)$.
\end{remark}

For the next lemma, recall the definition of $\qut^*_{c_0}(G')$ from Definition \ref{defQut*}. Note that a matching in a graph is a set of non-adjacent edges, i.e. no two edges share a common vertex. 

\begin{lemma}\label{lemcol}%-1 not in 1^{S_l \cap S_k}
Let $H$ be a graph and $G_*(M_H,0)$ as in Definition \ref{defnewgraph}. We choose the color of the added paths in $G_*'(M_H,0)$ to be $-1$, i.e. the same color as any edge between $V_l$ and $V_k$ for $k \neq l$. Choosing $c_0=-1$, we get $\qut^*_{-1}(G_*'(M_H,0))\cong \qut(G_*'(M_H,0))$.
\end{lemma}

\proof
We start by showing that every edge-color class except $c_0$ in $G_*'(M_H,0)$ is a matching. For this, first note that every edge between $V_k$ and $V_l$, $k\neq l$ as well as every edge in an added path has color $c_0$. Thus, the remaining edges are between vertices $(k,\alpha)$, $(k, \beta)$ in $V_k$ for some $k$. Thus, the color of the edge is the function $\alpha\triangle \beta: S_k\to \{\pm1\}^{S_k}$, $(\alpha\triangle\beta)_i=\alpha_i\beta_i$. Fix a vertex $(k,\alpha)$. For every $\beta\neq \beta'$ it holds $\alpha\triangle \beta \neq\alpha\triangle \beta'$, as otherwise $\beta_i=\beta_i'$ for all $i$. Thus, no two edges of the same color share a common vertex, which means that the color classes are matchings.

Let $c((i,k))= c((j,l))\neq c_0$. Then
\begin{align*}
    u_{ij}u_{kl}=u_{ij}\left(\sum_{a;c((j,a))=c((i,k))}u_{ka}\right)=u_{ij}\left(\sum_{a}u_{ka}\right)=u_{ij}, %=\left(\sum_{a;c((j,a))=c((i,k))}u_{ka}\right)u_{ij}=u_{kl}u_{ij},
\end{align*}
because $l$ is the only neighbor of $j$ such that $c((i,k))= c((j,l))$, since the color classes are matchings. Using the adjoint then yields $u_{ij}u_{kl}=u_{kl}u_{ij}$. We deduce $\qut^*_{-1}(G_*'(M_H,0))\cong \qut(G_*'(M_H,0))$. 
\qeds

We get the following theorem. 

\begin{theorem}\label{thm:qauteqsolgroup}
Let $H$ be a graph with $\deg(v)\geq 2$ for $v\in V(H)$ and $G_*(M_H,0)$ be as in Definition \ref{defnewgraph}. Let $G_*''(M_H,0)$ be the graph as in Lemma \ref{lemcol}. Then $\qut(G_*''(M_H,0))\cong \widehat{\Gamma_H}$. 
\end{theorem}

\proof
First note that we have $\deg(v)\geq 3$ for $v\in V(G_*(M_H,0))$ (since $\deg(w)\geq 2$ for $w\in V(H)$, we see that every $v\in V(G_*(M_H,0))$ has at least one neighbor that is associated to the same equation and at least two neighbors that come from different equations but share one of the at least two variables in the equation). We also know $\qut^*_{-1}(G_*'(M_H,0))\cong \qut(G_*'(M_H,0))$ by Lemma \ref{lemcol}. Thus, Corollary \ref{coreq} yields $\qut(G_*''(M_H,0))\cong \qut(G_*(M_H,0))$. Using Theorem \ref{thm:qaut2cgamma}, we deduce $\qut(G_*''(M_H,0))\cong \widehat{\Gamma_H}$.
\qed\vspace{0.5cm}

Using the theorem above, we are able to obtain a graph with quantum symmetry and finite quantum automorphism group. This graph is constructed from the constraint system associated to $K_{3,4}$.

\begin{cor}\label{cor:qsymfinqaut}
The graph $G_*''(M_{K_{3,4}},0)$ has quantum symmetry and finite quantum automorphism group. 
\end{cor}

\proof
It is easy to see that $\deg(v)\geq 2$ for all $v\in V(K_{3,4})$. Therefore $C(\qut(G_*''(M_{K_{3,4}},0))\cong C^*(\Gamma_{K_{3,4}}))$ by the previous theorem. Furthermore, it was shown in \cite{Paddock} (by using Sage \cite{sagemath}) that $\Gamma_{K_{3,4}}$ is finite and non-abelian, which yields that $C^*(\Gamma_{K_{3,4}})$ is finite-dimensional and non-commutative. By the isomorphism above, we get that $C(\qut(G_*''(M_{K_{3,4}},0)))$ is \linebreak finite-dimensional and non-commutative which yields the assertion.
\qed

\section{Quantum Isomorphisms}\label{sec:qiso}

 In this section we consider quantum isomorphisms between the colored graphs $G(M,b)$ and $G(M,b')$ for $b \ne b'$. In particular we give analogs of Lemma~\ref{lem:qutGdecomp} and Theorem~\ref{thm:qaut2cgamma} for this case. From this we are able to obtain non-isomorphic colored graphs $G(M,b)$ and $G(M,b')$ that are quantum isomorphic but neither $G(M,b)$ nor $G(M,b')$ has quantum symmetry. Further, applying decoloring techniques as in Section~\ref{sec:decolor} allows us to produce non-isomorphic uncolored graphs $G$ and $G'$ without quantum symmetry that are nevertheless quantum isomorphic. This appears to be the first known such example.

To define quantum isomorphism of colored graphs, we must provide a suitable generalization of the isomorphism game. The way to do this is quite natural.

\begin{definition}[Isomorphism game for colored graphs]
Given colored graphs $G$ and $G'$, with respective color functions $c$ and $c'$, the $(G,G')$-isomorphism game has as inputs and outputs for both Alice and Bob the set $V(G) \cup V(G')$\footnote{We are implicitly assuming here that the vertex sets of $G$ and $G'$ are disjoint. Alternatively we can assume that as part of their input they are told which graph their input should be associated to, and similarly for their output.}. To win, upon receiving $x \in V(G)$ Alice (respectively Bob) must respond with $y \in V(G')$ and vice versa. If this condition is met, then there is a vertex $g_a \in V(G)$ that is either Alice's input or output, and there is similarly $g_b \in V(G)$ and $g'_a, g'_b \in V(G')$. Alice and Bob then win if the following conditions are met:
\begin{enumerate}
    \item $c(g_a) = c'(g'_a)$ and $c(g_b) = c'(g'_b)$;
    \item $g_a = g_b$ if and only if $g'_a = g'_b$;
    \item $(g_a,g_b)$ is an edge of color $c$ if and only if $(g'_a,g'_b)$ is an edge of color $c$.
\end{enumerate}
\end{definition}

We then say that two colored graphs $G$ and $G'$ are \emph{quantum isomorphic} if there is a quantum strategy\footnote{There are a few different models of quantum strategies that can be considered. In this work we are always referring to the commuting operator framework when we mention quantum strategies.} that wins the $(G,G')$-isomorphism game with probability 1. In~\cite{qperms} it was shown that (uncolored) graphs $G$ and $G'$ are quantum isomorphic if and only if there exists a magic unitary $u$ such that $A_Gu = uA_{G'}$. Precisely the same proof applied to colored graphs gives the following:

\begin{prop}\label{prop:qisocolor}
Colored graphs $G$ and $G'$ with coloring functions $c$ and $c'$ are quantum isomorphic if and only if there exists a $V(G) \times V(G')$ magic unitary $u$ such that $c(g) \ne c'(g')$ implies $u_{gg'} = 0$ for $g \in V(G)$, $g' \in V(G')$, and $A_{G_c} u = uA_{G'_c}$ for all edge colors $c$.
\end{prop}

Based on the above, we define the following\footnote{This $C^*$-algebra is analogous to the $*$-algebra of a nonlocal game defined in~\cite{helton2019}, applied to the isomorphism game.}:

\begin{definition}[Isomorphism algebra]\label{def:isoalg}
Given colored graphs $G$ and $G'$ with coloring functions $c$ and $c'$, we define the \emph{isomorphism algebra}, denoted $Iso(G,G')$ as the universal $C^*$-algebra generated by elements $u_{gg'}$ for $g \in V(G)$, $g' \in V(G')$ such that
\begin{align}
    &u_{gg'}^2 = u_{gg'}^* = u_{gg'} \text{ for all } g \in V(G), g' \in V(G');\label{isoalg1} \\
    &\sum_{h \in V(G)} u_{hg'} = 1 = \sum_{h' \in V(G')} u_{gh'} \text{ for all } g \in V(G), g' \in V(G'); \label{isoalg2}\\
    &u_{gg'} = 0 \text{ if } c(g) \ne c'(g'); \label{isoalg3}\\
    &A_{G_c} u = u A_{G'_c} \text{ for all edge colors } c.\label{isoalg4}
\end{align}
\end{definition}

It is clear from Proposition~\ref{prop:qisocolor} and the above definition that the following holds:

\begin{prop}
Let $G$ and $G'$ be colored graphs. Then $G \cong_q G'$ if and only if $Iso(G,G')$ is nontrivial.
\end{prop}

We will also need to define the following algebra for our results in this section:

\begin{definition}\label{def:solgroupquotient}
For $M \in \mathbb{F}_2^{m \times n}$ and $b \in \mathbb{F}_2^m$, we define $\mathcal{A}(M,b)$ to be the universal $C^*$ algebra generated by $x_1, \ldots, x_n$ such that
\begin{enumerate}
    \item $x_i^2 = 1$ for all $i \in [n]$;
    \item $x_ix_j = x_jx_i$ if there exists $k \in [m]$ such that $M_{ki} = M_{kj} = 1$;
    \item $\prod_{i : M_{ki} = 1} x_i = (-1)^{b_k}$ for all $k \in [m]$.
\end{enumerate}
\end{definition}

Note that $\mathcal{A}(M,b)$ is not $C^*(\Gamma)$ for $\Gamma = \Gamma(M,b)$, unless $b = 0$. Rather, for $b \ne 0$, the algebra $\mathcal{A}(M,b)$ is isomorphic to the algebra $pC^*(\Gamma(M,b)) := \{px : x \in C^*(\Gamma(M,b))\}$ where $p = \frac{1}{2}(1-\gamma)$, which may be trivial for some $M$ and $b$. We will show that $Iso(G(M,b),G(M,b'))$ is isomorphic to $\mathcal{A}(M,b+b')$, but first we need an analog of Lemma~\ref{lem:qutGdecomp}. We omit the proof because it is almost identical to the proof of Lemma~\ref{lem:qutGdecomp}. The only difference is that one index of the magic unitary $u$ is a vertex of $G(M,b)$ and the other is a vertex of $G(M,b')$ instead of both being from the same graph.
     
% For uncolored graphs $G$ and $G'$, it is said that they are quantum isomorphic if there exists a quantum strategy\footnote{There are a few different models of quantum strategies that can be considered. In this work we are always referring to the commuting operator framework when we mention quantum strategies.} that wins the $(G,G')$-isomorphism game with probability 1. In~\cite{qperms} it was shown that this is equivalent to the existence of a magic unitary $u$ such that $A_G u = u A_{G'}$. For colored graphs, $G$ and $G'$, we can define the $(G,G')$-isomorphism game to have inputs $V(G) \cup V(G')$ and 

\begin{lemma}\label{lem:qisodecomp}
Let $M \in \mathbb{F}_2^{m \times n}$ and $b,b' \in \mathbb{F}_2^m$. Set $G = G(M,b)$ and $G' = G(M,b')$. If $u$ is a $V(G) \times V(G')$ magic unitary that gives a quantum isomorphism from $G$ to $G'$ then $u$ has the form
%\[u = \bigoplus_{k \in [m]} u^{(k)},\]
\begin{align}
\begin{blockarray}{cccccc}
&V_1(G')&V_2(G')&\dots&\dots&V_m(G')\\
\begin{block}{c(ccccc)}
V_1(G)&u^{(1)}&0&0&\dots&0\\ V_2(G)&0&u^{(2)}&0&\dots&0\\ \vdots&0&0&u^{(3)}&\dots&0\\\vdots&\vdots &\vdots&\vdots&\ddots&0\\V_m(G)&0&0&0&0&u^{(m)} \\
\end{block}
\end{blockarray}
\end{align}
where each $u^{(k)}$ is a magic unitary. Furthermore, $u^{(k)}_{\alpha,\beta} := u_{(k,\alpha),(k,\beta)}$ depends only on $k$ and the value of $\alpha \triangle \beta$ and therefore the entries of $u^{(k)}$ pairwise commute for each $k \in [m]$.
\end{lemma}

We are now able to prove the following partial analog of Theorem~\ref{thm:qaut2cgamma}:
\begin{theorem}\label{thm:qiso2cgamma}
Let $M \in \mathbb{F}_2^{m \times n}$ and $b \ne b' \in \mathbb{F}_2^m$. Set $G = G(M,b)$, $G' = G(M,b')$, and $\mathcal{A} = \mathcal{A}(M,b+b')$. Then there exists an isomorphism $\varphi : \mathcal{A} \to Iso(G,G')$. 
\end{theorem}
\proof
The proof is very similar to that of Theorem~\ref{thm:qaut2cgamma}. We will first give a homomorphism $\varphi_1$ from $\mathcal{A}$ to $Iso(G,G')$ and then a homomorphism $\varphi_2$ in the other direction, and finally show that these are inverses of each other and therefore isomorphisms.\\

%First note that $p \in C^*(\Gamma)$ is a projection that commutes with every element of $C^*(\Gamma)$, the latter following from the fact that $\gamma$ commutes with every element of $\Gamma$. It follows that $p$ is the identity in $\mathcal{A}$.\\

\noindent\emph{Step 1: Construction of a $*$-homomorphism $\varphi_1: \mathcal{A} \to Iso(G,G')$.}\\
For this step, we will construct elements $y_i$ of $Iso(G,G')$ that satisfy the relations of the generators $x_i$ of $\mathcal{A}$. This proves that there is a $*$-homomorphism $\varphi_1$ from $\mathcal{A}$ to $Iso(G,G')$ such that $\varphi_1(x_i) = y_i$. Later we will see that $\varphi_1$ is in fact an isomorphism.

Let $u$ be magic unitary of generators for $Iso(G,G')$. By Lemma~\ref{lem:qisodecomp}, $u = \bigoplus_{k \in [m]} u^{(k)}$ where each $u^{(k)}$ is a $V_k(G) \times V_k(G')$ magic unitary. Moreover, the operator $u^{(k)}_{\alpha,\beta} := u_{(k,\alpha),(k,\beta)}$ depends only on $k$ and the value of $\alpha \triangle \beta \in \pm 1^{S_k}_{b_k+b_k'}$. Thus, as in the proof of Lemma~\ref{lem:qutGdecomp}, for each $\delta \in \pm 1^{S_k}_{b_k+b'_k}$ we let $u^{(k)}_\delta$ denote $u^{(k)}_{\alpha,\beta}$ such that $\alpha \triangle \beta = \delta$. Note that for any $\alpha \in \pm 1^{S_k}_{b_k}$, we have that
\begin{equation}\label{eq:udeltares2}
    \{u^{(k)}_\delta : \delta \in \pm 1^{S_k}_{b_k+b'_k}\}  = \{u^{(k)}_{\alpha,\beta} : \beta \in \pm 1^{S_k}_{b'_k}\},
\end{equation}
and thus every row (and similarly every column) of $u^{(k)}$ contains the same set of operators.

Now we define $y^{(k)}_i = \sum_{\alpha \in \pm 1^{S_k}_{b_k+b'_k}} \alpha_i u^{(k)}_\alpha$ for all $k \in [m]$ and $i \in S_k$. We first aim to show that $y^{(k)}_i$ does not depend on $k$.

Fix $i \in [n]$ and $l,k \in [m]$ such that $i \in S_l \cap S_k$. For each $\delta \in \pm 1^{S_l \cap S_k}$, let $A_\delta$ (resp.~$A'_\delta$) be the adjacency matrices of the graph consisting of the edges of $G$ (resp.~$G'$) colored $\delta$. Further, let $B_\delta$ (resp. $B'_\delta$) be the submatrix of $A_\delta$ (resp.~$A'_\delta$) consisting of the rows indexed by $V_l(G)$ and columns indexed by $V_k(G)$ (resp.~$V_l(G')$ and $V_k(G')$). In other words,
\[(B_\delta)_{(l,\alpha),(k,\beta)} = \begin{cases}1 & \text{if } \alpha \triangle \beta = \delta \\ 0 & \text{o.w.}\end{cases},\]
and similarly for $B'_\delta$. Also define $A_+ = \sum_{\delta \in \pm 1^{S_l \cap S_k}, \delta_i = +1} A_\delta$ and $A_- = \sum_{\delta \in \pm 1^{S_l \cap S_k}, \delta_i = -1} A_\delta$, and similarly define $A'_+$, $A'_-$, $B_+$, $B_-$, $B'_+$, and $B'_-$. Then
\[(B_+)_{(l,\alpha),(k,\beta)} = \begin{cases}1 & \text{if } \alpha_i = \beta_i \\ 0 & \text{o.w.}\end{cases} \quad \text{ \& } \quad (B_-)_{(l,\alpha),(k,\beta)} = \begin{cases}1 & \text{if } \alpha_i = -\beta_i \\ 0 & \text{o.w.}\end{cases}\]
and similarly for $B'_+$ and $B'_-$. Since $u$ must satisfy $A_\delta u = uA'_\delta$ by definition, it must also satisfy $A_\pm u = uA'_\pm$. This implies that $B_+ u^{(k)} =  u^{(l)}B'_+$ and $B_- u^{(k)} = u^{(l)}B'_-$. Considering the $(l,\alpha),(k,\beta)$ entry of both sides of the former equation, we see that
\[\sum_{\substack{\beta' \in \pm 1^{S_k}_{b_k} \\ \text{s.t. } \beta'_i =  \alpha_i}}u^{(k)}_{\beta',\beta} = \sum_{\substack{\alpha' \in \pm 1^{S_l}_{b'_l} \\ \text{s.t. } \alpha'_i = \beta_i}}u^{(l)}_{\alpha,\alpha'}.\]
Note that for every term in the above sums, we have that $\alpha_i\alpha'_i = \alpha_i\beta_i = \beta'_i\beta_i$. Thus, reexpressing $u^{(l)}_{\alpha,\alpha'}$ as $u^{(l)}_{\tilde{\alpha}}$ where $\tilde{\alpha} = \alpha \triangle \alpha'$, and similarly for $u^{(k)}_{\beta',\beta}$ we have that
\[\sum_{\substack{\tilde{\beta} \in \pm 1^{S_k}_{b_k+b'_k} \\ \text{s.t. } \tilde{\beta}_i =  \alpha_i\beta_i}}u^{(k)}_{\tilde{\beta}} = \sum_{\substack{\tilde{\alpha} \in \pm 1^{S_l}_{b_l+b'_l} \\ \text{s.t. } \tilde{\alpha}_i = \alpha_i\beta_i}}u^{(l)}_{\tilde{\alpha}}.\]
Doing the same for $B_- u^{(k)}= u^{(l)}B'_-$ yields the same equation but with $\alpha_i\beta_i$ replaced with $-\alpha_i\beta_i$ and combining these proves that
\[\sum_{\tilde{\beta} \in \pm 1^{S_k}_{b_k + b'_k}, \tilde{\beta}_i = x}u^{(k)}_{\tilde{\beta}} = \sum_{\tilde{\alpha} \in \pm 1^{S_l}_{b_l + b'_l}, \tilde{\alpha}_i = x}u^{(l)}_{\tilde{\alpha}}\]
for any $x \in \{+1,-1\}$. Therefore,
\[y^{(l)}_i = \left(\sum_{\substack{\alpha \in \pm 1^{S_l}_{b_l + b'_l} \\ \text{s.t. } \alpha_i = +1}} u^{(l)}_{\alpha}\right) - \left(\sum_{\substack{\alpha \in \pm 1^{S_l}_{b_l + b'_l} \\ \text{s.t. } \alpha_i = -1}} u^{(l)}_{\alpha}\right) = \left(\sum_{\substack{\alpha \in \pm 1^{S_k}_{b_k + b'_k} \\ \text{s.t. } \alpha_i = +1}} u^{(k)}_{\alpha}\right) - \left(\sum_{\substack{\alpha \in \pm 1^{S_k}_{b_k + b'_k} \\ \text{s.t. } \alpha_i = -1}} u^{(k)}_{\alpha}\right) = y^{(k)}_i.\]
So we have shown that the value of $y^{(k)}_i$ does not depend on $k$, and thus we will simply denote this operator by $y_i$.

Now note that since $y_i$ is a linear combination (with real coefficients) of the operators $u^{(k)}_\alpha$ for $k \in [m]$ such that $i \in S_k$, and these operators are entries of the magic unitary $u$, we have that $y_i^* = y_i$. Also, by~\eqref{eq:udeltares2} we have that $u^{(k)}_\alpha u^{(k)}_\beta = 0$ for $\alpha \ne \beta$ and thus
\[y_i^2 = \left(\sum_{\alpha \in \pm 1^{S_k}_0} \alpha_i u^{(k)}_\alpha\right)^2 = \sum_{\alpha \in \pm 1^{S_k}_0} \alpha^2_i \left(u^{(k)}_\alpha\right)^2 = \sum_{\alpha \in \pm 1^{S_k}_0} u^{(k)}_\alpha = 1.\]
Thus the $y_i$ satisfy relation (1) from Definition~\ref{def:solgroupquotient}.

Next we will show that relation (2) of Definition~\ref{def:solgroupquotient} holds, i.e., that $y_iy_j = y_jy_i$ if there exists $k \in [m]$ such that $i,j \in S_k$. Suppose that $i,j,k$ are as described. Then $y_i = y_i^{(k)}$ and $y_j = y_j^{(k)}$ are both linear combinations of the entries of $u^{(k)}$ which pairwise commute by Lemma~\ref{lem:qisodecomp}. Therefore $y_i$ and $y_j$ commute as desired.

Lastly, we must show that relation (3) of Definition~\ref{def:solgroupquotient} holds, i.e., that $\prod_{i \in S_k}y_i = (-1)^{b_k + b'_k}$ for all $k \in [m]$. We have that
\begin{align*}
    \prod_{i \in S_k}y_i &= \prod_{i \in S_k} y^{(k)}_i \\
    &= \prod_{i \in S_k} \left(\sum_{\alpha \in \pm 1^{S_k}_{b_k + b'_k}} \alpha_i u^{(k)}_\alpha\right)
\end{align*}
Since the elements $u^{(k)}_\alpha$ are pairwise orthogonal, when we expand the above product all cross terms disappear and we obtain
\[\prod_{i \in S_k}y_i = \sum_{\alpha \in \pm 1^{S_k}_{b_k + b'_k}} \left(\prod_{i \in S_k} \alpha_i\right)u^{(k)}_\alpha = \sum_{\alpha \in \pm 1^{S_k}_{b_k + b'_k}} (-1)^{b_k + b'_k} u^{(k)}_\alpha = (-1)^{b_k + b'_k},\]
since $\prod_{i \in S_k} \alpha_i = (-1)^{b_k + b'_k}$ for all $\alpha \in \pm 1^{S_k}_{b_k + b'_k}$ by definition. Therefore the elements $y_i \in Iso(G,G')$ for $i \in [n]$ satisfy the relations of the generators of $\mathcal{A}$ and thus there exists a $*$-homomorphism $\varphi_1$ from $\mathcal{A}$ to $Iso(G,G')$ such that $\varphi_1(x_i) = y_i$ for all $i \in [n]$.\\

\noindent\emph{Step 2: Construction of a $*$-homomorphism $\varphi_2: Iso(G,G') \to \mathcal{A}$.}\\
This step is almost identical to Step 2 of the proof of Theorem~\ref{thm:qaut2cgamma}, and so we omit it. We only remark that the biggest change is that when showing that $v_\alpha^{(k)} = 0$ for $\alpha \notin \pm 1^{S_k}_{b_k + b'_k}$, the analog of Equation~\eqref{eq:vka0} gets an additional factor of $(-1)^{b_k + b'_k}$ in the first and last expressions.\\

\noindent\emph{Step 3: Showing that $\varphi_1$ and $\varphi_2$ are inverses of each other.}\\
This step is nearly identical to Step 3 of the proof of Theorem~\ref{thm:qaut2cgamma} and so we omit it.\qeds

\begin{cor}\label{cor:qiso2solgroup}
Let $M \in \mathbb{F}_2^{m \times n}$ and $b \ne b' \in \mathbb{F}_2^m$. Then the following are equivalent:
\begin{enumerate}
    \item there is a quantum isomorphism from $G(M,b)$ to $G(M,b')$;
    \item the algebra $\mathcal{A}(M,b+b')$ is nontrivial;
    \item the element $\gamma$ of the solution group $\Gamma(M,b+b')$ is not equal to the identity.
\end{enumerate}
\end{cor}
\proof
Proposition~\ref{prop:qisocolor} states that (1) is equivalent to $Iso(G(M,b),G(M,b'))$ being nontrivial, which is equivalent to (2) by Theorem~\ref{thm:qiso2cgamma}. To see that (2) and (3) are equivalent, recall that $\mathcal{A}(M,b+b')$ is isomorphic to $pC^*(\Gamma(M,b+b')) = \{px : x \in nC^*(\Gamma(M,b+b'))\}$ where $p = \frac{1}{2}(1-\gamma)$, and so these algebras are nontrivial if and only if $\gamma \ne 1$ in $C^*(\Gamma(M,b+b'))$ which is equivalent to $\gamma$ not being equal to the identity in $\Gamma(M,b+b')$.\qeds

Note that the equivalence of $(1)$ and $(3)$ in the corollary above differs from earlier results (e.g., from~\cite{nonsignalling}) because the graphs $G(M,b)$ and $G(M,b')$ here are \emph{colored} graphs, rather than uncolored ones.

We are now able to give two non-isomorphic colored graphs without quantum symmetry that are quantum isomorphic. Recall the definition of $G_*(M,0)$ from Definition \ref{defnewgraph}.%comment why we choose G_* (otherwise they could have q sym)

\begin{cor}\label{cor:qisonoqsym}
Let $M$ be the incidence matrix of $K_{3,3}$, and let $b \in \mathbb{F}_2^{V(K_{3,3})}$ be any vector with exactly one entry equal to 1. Then $G_*(M,0)$ and $G_*(M,b)$ are non-isomorphic colored graphs with no quantum symmetry which are nevertheless quantum isomorphic.
\end{cor}
\proof
In this case it is well known that the element $\gamma$ of the solution group $\Gamma(M,b+0) = \Gamma(M,b)$ is not equal to the identity~\cite{CM}, since this is the solution group of the Mermin-Peres magic square game. It then follows from Corollary~\ref{cor:qiso2solgroup} that $G_*(M,0)$ and $G_*(M,b)$ are quantum isomorphic. On the other hand, it is known that even the uncolored versions of these graphs are non-isomorphic~\cite{nonsignalling}, and thus the colored versions are as well. Lastly, it was shown in~\cite{Paddock} that $\Gamma(M,0)$ is abelian, and therefore by Theorem~\ref{thm:qaut2cgamma} the algebras $C(\qut(G_*(M,0)))$ and $C(\qut(G_*(M,b)))$ are commutative, i.e., the graphs $G_*(M,0)$ and $G_*(M,b)$ have no quantum symmetry.\qeds

\begin{cor}\label{cor:qisonoqsymnocolor}
Let $M$ be the incidence matrix of $K_{3,3}$, and let $b \in \mathbb{F}_2^{V(K_{3,3})}$ be any vector with exactly one entry equal to 1. Let $G_*''(M,0)$ and $G_*''(M,b)$ the graphs constructed from $G_*(M,0)$ and $G_*(M,b)$ as in Lemma \ref{lemcol}, where we choose the added paths in both graphs to have the same length for vertices of the same color and edges of the same color. Then $G_*''(M,0)$ and $G_*''(M,b)$ are non-isomorphic graphs with no quantum symmetry which are nevertheless quantum isomorphic. 
\end{cor}

\proof
We start by showing that $G_*''(M,0)$ and $G_*''(M,b)$ have no quantum symmetry. By the previous corollary, we know that $G_*(M,0)$ and $G_*(M,b)$ have no quantum symmetry. We have $C(\qut^*_{-1}(G_*'(M,0)))\cong C(\qut(G_*'(M,0)))$ and $C(\qut^*_{-1}(G_*'(M,b)))\cong C(\qut(G_*'(M,b)))$ by Lemma \ref{lemcol}. Therefore, Corollary \ref{coreq} yields that $G_*''(M,0)$ and $G_*''(M,b)$ have no quantum symmetry. 

It remains to show that $G_*''(M,0)$ and $G_*''(M,b)$ are quantum isomorphic. By Corollary \ref{cor:qisonoqsym}, we know that $G_*(M,0)$ and $G_*(M,b)$ are quantum isomorphic. Thus, there exists a $C^*$-algebra $A$ and a matrix $U\in M_n(A)$ such that Relations \eqref{isoalg1}--\eqref{isoalg4} hold. In $G_*''(M,0)$, we denote the vertex that subdivided the edge $e$ in $G_*(M,0)$ by $e_0$ and the vertices in the added path with $d(e_0,e_i)=i$ by $e_i$, $1 \le i \le m_{c(e)}$. Moreover, we denote by $v_i$, $0 \le i \le n_{c(v)}$ the vertices in the added path with $d(v,v_i)=i$. We similarly label the vertices of $G_*''(M,b)$. 

Note that by a similar reasoning as in Lemma \ref{lemcol}, we get $u_{ij}u_{kl}=u_{kl}u_{ij}$ for all $e=(i,k), f = (j,l)$ with $c(e)\neq c_0 \neq c'(f)$. Define the magic unitary $W$ as follows: 
\begin{align*}
    w_{i_kj_l}&= \delta_{kl}u_{ij} \text{ for }i \in V(G(M,0)), j\in V(G(M,b)),\\ %0\le k \le n_{c(i)}, 0\le l \le n_{c(j)},
    w_{e_kf_l}&=\delta_{kl}(u_{ac}u_{bd}+u_{bd}u_{ac}) \text{ for }e=(a,b)\in E(G(M,0)), f=(c,d)\in E(G(M,b)),\\
    w_{i_kf_l}&=0 \text{ for }i \in V(G(M,0)), f\in E(G(M,b)),\\
    w_{e_kj_l}&=0 \text{ for }e \in E(G(M,0)), j\in V(G(M,b)), 
\end{align*}
for all occurring $k,l$. Similar to the proofs of Theorem \ref{thmvertexcolor} and Theorem \ref{theoremedgecolor}, we see that $A_{G_*''(M,0)}W=WA_{G_*''(M,b)}$. We conclude that $G_*''(M,0)$ and $G_*''(M,b)$ are quantum isomorphic. %more detailed? Only A_{G_c}(4.5) and A_{G_0}(4.8) differ, where we only use the commutation with u for both.
\qeds

\paragraph{Acknowledgments.}\phantom{a}\newline
\noindent\textbf{DR} was partially supported by the VILLUM Centre of Excellence for the 
Mathematics of Quantum Theory (QMATH) project no. 35544.\\
\noindent\textbf{SSch} would like to thank Matthew Daws and Christian Voigt for helpful discussions. This project has received funding from the European Union's Horizon 2020 research and innovation programme under the Marie Sklodowska-Curie grant agreement No. 101030346.

\bibliographystyle{plainurl}
\bibliography{graphqsymfin}

\end{document}